\documentclass[11pt]{article}
\usepackage{amsfonts,amsmath,amssymb,amsthm,amscd,epsfig,epic,eepic}
\usepackage{graphics}
\usepackage{fullpage}
\usepackage{pstricks}

\usepackage[latin1]{inputenc}

\usepackage{graphicx,color}
\usepackage{wrapfig}
\usepackage{multirow}
\usepackage{array}
\usepackage{multicol}
\usepackage{amsmath,amssymb}
\usepackage{float}
\usepackage{mathrsfs}

\usepackage[all]{xy}
\xyoption{matrix}
\xyoption{arrow}

\def\ds{\displaystyle}

\def\cp1{\projectif^1}
\def\riemann{\widehat{\complexes}}

\def\eps{\varepsilon}
\def\ph{\varphi}

\def\rah{\begin{eqnarray}}
\def\arh{\end{eqnarray}}

\def\rae{\begin{eqnarray*}}
\def\are{\end{eqnarray*}}

\def\demi{\frac{1}{2}}
\def\fr{\frac}

\def\tend{\rightarrow}
\def\sans{\backslash}

\def\commutatif{\ar@{}[rd]|{\circlearrowleft}}

\newcommand{\barre}[1]{\overline{#1}}
\newcommand{\inv}[1]{\frac{1}{#1}}

\DeclareMathAlphabet{\mathpzc}{OT1}{pzc}{m}{it}
\renewcommand{\Re}{\operatorname{Re}}
\renewcommand{\Im}{\operatorname{Im}}

\def\parel{\Re} 
\def\parim{\Im}

\newtheorem{thm}{\bf Theorem}[section]
\newtheorem{lem}[thm]{\bf Lemma}

\newtheorem{prop}[thm]{\bf Proposition}

\newtheorem{cor}[thm]{\bf Corollary}

\newtheorem{rem}[thm]{\bf Remark}
\newtheorem{conjecture}[thm]{\bf Conjecture}

\def\nin{\notin}

\def\preuve{\ \\\textit{Proof. }}
\def\cqfd{\hfill $\square$}

\def\entiers{\mathbb{N}}
\def\naturels{\entiers}
\def\relatifs{\mathbb{Z}}

\def\reels{\mathbb{R}}
\def\complexes{\mathbb{C}}

\def\projectif{\mathbb{P}}

\def\ii{i}
\def\disque{\mathbb{D}}

\def\mandelbrot{\mathcal{M}}

\def\bd{\barre{\disque}}

 \evensidemargin=\oddsidemargin

\makeatletter
\gdef\thmhead@plain#1#2#3{%
  \thmname{#1}\thmnumber{\@ifnotempty{#1}{ }#2}%
  \thmnote{ {\mdseries#3}}}
\let\thmhead\thmhead@plain
\makeatother
\theoremstyle{plain}
\newtheorem{theorem}{Theorem}[section]

\newtheorem{question}[theorem]{Question}

\newtheorem{pro}{Proposition}

\theoremstyle{definition}
\newtheorem{definition}[theorem]{Definition}

\theoremstyle{remark}

\def\alinea#1{\hfill\break%
  \hbox to \parindent{\hss{\upshape{\bf #1)}}\enspace}\ignorespaces}
\def\bul{\hfill\break\hbox to\parindent{\hss$\bullet$\enspace}\ignorespaces}

\def\id{\hbox{Id}}

\def\C{\mathbb{C}}
\def\D{\mathbb{D}}

\def\M{\mathcal{M}}

\def\N{\mathbb{N}}
\def\Q{\mathbb{Q}}
\def\R{\mathbb{R}}
\def\S{\mathbb{S}}

\def\Z{\mathbb{Z}}

\def\Cap_#1{\bigcap\limits_{#1}}
\def\Cup_#1{\bigcup\limits_{#1}}
\def\ol{\overline}
\def\ul{\underline}

\makeatletter
\def\cqfdsymb{\relax\protect\ifmmode\else\unskip\nobreak\fi
\quad\hfill$\bgroup
\vcenter{\hrule\hbox{\vrule\@height.6em\kern.6em\vrule}\hrule}\egroup$}
\def\cqfd{\cqfdsymb  \endtrivlist}
\gdef\rom#1{\leavevmode\skip@\lastskip\unskip\/%
        \ifdim\skip@=\z@\else\hskip\skip@\fi{\normalshape#1}}
\makeatother

\def\ueps{{\ul{\eps}}}

\begin{document}
\title{{\bf On (non-)local-connectivity of some Julia sets}}

\author{{\sc A. Dezotti}, {\sc P. Roesch}\thanks{Institut of Mathematics of Toulouse}}

\date{\today}
\maketitle
\begin{abstract} This article deals with the question of local connectivity of the Julia set of 
polynomials and rational maps. It essentially presents conjectures and questions. 
\end{abstract}

\section*{Introduction}

In this note we discuss  the following question : When is the Julia set of a rational map connected but not locally-connected?
We propose some  conjectures and develop a  model of non-locally-connected Julia sets in the case of infinitely renormalizable quadratic polynomials, a situation where one hopes to find a precise answer.

The question of local connectivity of the Julia set has been studied extensively for quadratic polynomials,
  but there is still no complete characterization of when a quadratic  polynomial has a
 connected and locally-connected Julia set. 
 In degree $2$, the question reduces  to the precise cases where the  polynomial has a Siegel disk or is infinitely renormalizable.
 J. Milnor proposed in his lecture~\cite{MiNlc}
 a quantitative condition to get a non-locally-connected Julia set which is infinitely satellite renormalizable.
 It follows  the work of  A. Douady and D. E. K. S{\o}rensen\,: In \cite{So}
  a description of the topological nature of a non-locally-connected Julia set is given, and some examples in the
 infinitely satellite renormalizable case are obtained.
 Nevertheless, the argument in \cite{So} is by continuity and 
gives no explicit condition.
G. Levin gave afterwards such a condition in \cite{Levin-2009} (see also Theorem\ref{thm de levin}).
 In section \ref{s:infinitely} of  this note, we present a model  of what the structure of the post-critical set in that setting  should be
 (it was originally created by X. Buff).

In a previous work, we  considered  polynomials of higher degrees. Here we present an example (section~\ref{s:example}) where  the local connectivity can be deduced by renormalization. 

The situation is even more  complicated for rational maps.
 Indeed, there are examples of rational maps  with Cremer points  such that the Julia set is locally-connected \cite{Ro4}.
 It seems more difficult to find examples of non-locally-connected Julia sets  in the space of rational maps.
 Nevertheless they exist, and can be easily obtained  by ``tuning'' 
from polynomials.  
From the way those rational maps are constructed,
 the natural  question  appears to be how much a rational map has to be related to a polynomial so that its Julia set is not locally-connected?
 Are polynomials pathological rational maps? 
In all cases presented here, when the rational map or polynomial has a connected but not locally-connected Julia set,
  a criterion is verified. We call it ``Douady-Sullivan criterion'' since it has been used the first time by them.

\vskip 0.51em 
{\it Aknowledgment} :
 The authors would like to thank the referee for his thorough reading of the manuscript and his comments that helped  improve the manuscript.
 Discussions with X. Buff, A. Ch\'eritat and Y. Yin, as well as remarks of  D. Cheraghi and  H. Inou inspired 
parts of this paper.

\section{Local connectivity }
\subsection{Generalities and first questions for polynomials}
Recall that the \textit{Julia set} of a rational map 
$f$ is the minimal totally invariant (under $f$ and $f^{-1}$) compact set containing at least $3$ points.
 Its  complement in the Riemann sphere,
 called the \textit{Fatou set}, is an open set whose components are all eventually periodic by Sullivan's Theorem.
 When the Julia set is connected, these components are all topological disks.
 Inside each of the periodic components the return map is conjugate near the boundary to some simple model (see~\cite{aut}).
If the boundary of the component is locally-connected, the model extends to the boundary (by Carath\'eodory's Theorem).
 One of the main reasons to consider the question of local connectivity for Julia sets is to get the model on the boundary.

 Recall that the models are given by the following maps from the unit disc to itself :\begin{itemize}
 \item 
 $Z_d(z)= z^d$, the {\it attracting case},  
 \item
 $B_d(z)=\frac{z^d+v}{1+v z^d}$  where  $v=\frac{d-1}{d+1}$, the  {\it parabolic case}, 
 \item 
$R_\theta(z)=e^{2i\pi \theta} z$, the  {\it Siegel case} (the corresponding  Fatou component is then called a {\it Siegel disk}).
   \end{itemize}
   In what follows we will always assume that the Julia sets considered are connected, even if it is not explicitly mentioned.
   
   \begin{lem}\label{l:decomposition}{\rm\cite[Th 4.4]{W}}
 The Julia set of a rational map is locally-connected if and only if the boundary of each Fatou component is locally-connected and
 for any $\epsilon>0$ only finitely many Fatou components have diameter greater than $\epsilon$.
 \end{lem}
 Hence, the question of whether each Fatou component has a locally-connected boundary is fundamental for a rational map.
 For a polynomial, the boundary of the unbounded Fatou component is the whole Julia set.
 Nevertheless, it is an interesting question to know if one can deduce some result looking only at the bounded Fatou components and their size.
 We will now give an answer  to this question. 
 
For polynomials let us recall the following result~\cite{RY} :
\begin{theorem}[(R-Yin)] \label{RY}Any bounded periodic Fatou component of a polynomial containing a critical point is a Jordan domain.
 \end{theorem} The following ``classical'' conjecture is the natural extension of this result to any bounded Fatou component. It has been proved recently in many cases by M. Shishikura.

 \begin{conjecture} The boundary of a periodic Siegel disk of a polynomial is always a Jordan curve.
 \end{conjecture}

A periodic point in the Julia set is called a {\it Cremer point} if the derivative of the return map at the fixed point is $e^{2i\pi t}$ with $t\in \R\setminus \Q$.
 Let us recall the following result (a proof will be sketched  in section~\ref{s:CremerSiegel})
\begin{pro}If a polynomial has either a Cremer periodic point or a periodic Siegel  disk with no critical point on the boundary (of the  cycle generated by the disk), then its Julia set is not locally-connected.
\end{pro}

This answers the question trivially since in the Cremer case there are no bounded Fatou components.
 One can also construct non-locally-connected Julia sets with  Fatou components 
 that are Jordan domains whose diameter tends to zero.
 Indeed, it is enough to take a polynomial containing both an attracting cycle and a Cremer point such that the orbit of the critical points does not accumulate on the boundary of the attracting basin.
 Then the attracting basin and all its pre-images are Jordan domains.
Moreover, using the  ``shrinking Lemma'' (see~\cite{TY} Prop. A.3 or \cite{LM} section $11.1$) in the complement of the post-critical set 
 ($\ol{P_f}$ where $P_f:=\ds{\bigcup_{c\in crit}\bigcup_{n\ge 1}}f^n(c)$)
 it is easy to see that the diameter of these Jordan domains goes to  zero.
 Such examples are easy to find in cubic families with one attracting fixed point (see~\cite{Ro3}).
 There we can find copies of the Mandelbrot set (the connectedness locus  for the quadratic family) in which we can choose  a doubly renormalizable restriction containing a  Cremer point (we will define renormalizable below in Definition~\ref{d:renorm}).

In the light of such examples, the previous question appears to be a naive one but its original motivation leads to the following  less naive question :

\begin{question} Let $P$ be a quadratic polynomial having a Siegel disk whose boundary is a Jordan curve containing the critical point. Is the Julia set locally-connected? 
\end{question}
To our knowledge there is no known counter-example, and in higher degrees the situation is even more complicated.
 One could imagine  to build a cubic polynomial  from a quadratic one that is non-locally-connected the following way :
One would have the Julia set of a quadratic polynomial  with a Siegel disk without the critical point on its boundary sitting in the Julia set of a cubic  polynomial.
 The critical point has to lie on some hairs around the Siegel disk.
 The idea is then to deform the cubic polynomial in the space of cubic polynomials in order to put the other critical point on the boundary of the Siegel disk. 
 This  kind of map has been  considered
 when the critical point belongs  to strict pre-images of the Siegel disk in~\cite{BuH}.
  One would get in our case a   polynomial with a Siegel disk containing one critical point on its boundary and another critical point on some hairs stemming from the Siegel disk. Nevertheless, it is not clear 
that having one critical point on the boundary of the Siegel disk will not force the hairs to disappear.
 Indeed, in the light of Douady-Sullivan criterion (see section~\ref{s:conj}),
 the non-local-connectivity seems to appear when the map presents some  injectivity,
 but here  around the boundary of  the Siegel disk the map is no longer injective, so there is no reason to expect the boundary to be topologically wild.

\begin{question}\label{q1}
Does there exist a non-renormalizable polynomial,  of degree $d\ge 3$, with a Siegel disk containing at least a critical point on its boundary and whose Julia set is not locally-connected?
\end{question}

Recall the definition of renormalizable maps.
 \begin{definition} \label{d:renorm}
A polynomial is said to be {\it renormalizable} if some iterate admits a polynomial-like restriction whose filled Julia set is connected.
 A {\it polynomial-like map} is a proper holomorphic map $f:U\to V$ where $U,V$ are topological disks with $\ol U\subset V$;
 one defines its filled Julia set as  $\bigcap f^{-n} (\ol U)$. (A polynomial is an example of polynomial-like map).
\end{definition}
Notice  that by  the following  connectedness principle (see \cite{mmullen-renorm}),
 in order for a renormalizable polynomial to have a locally-connected Julia set, the Julia set of its renormalized map should be locally-connected.
\begin{theorem}(Connectedness principle)\label{cp}
Let $f:\C\to \C$  be a polynomial with connected filled Julia set $K(f)$.  Let $f^n: U \to V$ be a renormalization of $f$ with filled  Julia set $K_n$. Then $\partial K_n\subset J(f)$ and 
for any closed connected set $L\subset K(f)$, $L\cap K_n$ is also connected. 
\end{theorem}

 Yoccoz proved that a quadratic polynomial that is finitely renormalizable and has only repelling periodic points has a locally-connected Julia set.
 In higher degree, the following result appears to be the most general known result. 
   
 \begin{theorem}[(Koslovski-van Strien)]\label{KvS}
The Julia set of a non-renormalizable polynomial without indifferent periodic points is locally-connected provided it is connected.
 \end{theorem}
One wonders if there is a way to combine previous results in order to justify the following : 
\begin{question}\label{q2}Is the Julia set of a polynomial  locally-connected provided that it is not infinitely renormalizable, it has no Cremer points and there is a critical point on the boundary of any cycle of Siegel disks?
\end{question}   
Question~\ref{q1} justifies partially question~\ref{q2}. Here are some further justifications.
 Using the work done in~\cite{PR}, one can construct  a puzzle in the  basins of parabolic cycles. 
Therefore, the proof of~\ref{KvS} (see~\cite{KV}) will adapt to the case of parabolic cycles as soon as one can construct with this puzzle a ``box mapping''.
 This can fail when the map is ``parabolic-like'' (see~\cite{Luna}).
 Nevertheless, the recent work of L. Lomonaco (\cite{Luna}) on parabolic-like maps should take care of this case.

In the finitely-remormalizable case, the proof proceeds by induction and uses the homeomorphism given by the straightening Theorem (\cite{DH2}).
 Theorem \ref{cubiclc} presented in next section is an example of this method.

Finally, let us point out that the case of infinitely renormalizable polynomials is much more subtle.
 Kahn, Levin, Lyubich, McMullen, van Strien ... gave conditions to obtain infinitely 
 renormalizable quadratic  polynomials with locally-connected Julia sets. Douady and S{\o}rensen gave examples of non-locally-connected infinitely renormalizable quadratic polynomials. We will discuss infinitely renormalizable polynomials in section~\ref{s:infinitelyDS} and in section~\ref{s:infinitely}.

 \subsection{From  Fatou components to the whole Julia set : an example}\label{s:example} 

We would like to end this section by a concrete example.
 We prove local connectivity of the Julia set of a polynomial knowing that it is renormalizable and that the small Julia set is locally-connected\,:
 
\begin{theorem}\label{cubiclc} Let $f_a(z)=z^{d-1}(z+da/(d-1))$ with $d\ge 3$ and  $a\in \C$,  be  the family of polynomials of degree $d$, with one fixed critical point of maximal multiplicity (up to affine conjugacy). Assume that the Julia set $J(f_a)$ is connected. If $f_a$ is renormalizable of lowest period $k$ around the ``free'' critical point $-a$,  we will denote by $Q_c$ the unique quadratic polynomial to which the restriction of $f_a^k$  is conjugate.  Then the Julia set $J(f_a)$ is locally-connected if and only if either $f_a$ is not renormalizable or $J(Q_c)$ is locally-connected.  
\end{theorem}
\proof First recall  that if $f_a$ or equivalently $Q_c$ is geometrically finite ({\it i.e.}  if the post-critical set intersects the Julia set at finitely many points), then the Julia sets are locally-connected by the result of~\cite{TY}.   We recall  in the following  the construction  of a graph ``adapted'' to the dynamics of the map $f_a$ (as presented in~\cite{Ro1}). Given a graph $\Gamma$  the connected component of $\C\setminus f^{-n}(\Gamma)$ containing $x$, is called the puzzle piece of depth $n$ containing $x$ and is denoted by  $P_n(x)$.   

\vskip 0.5em 
$*$ {\it Claim : There exists a graph $\Gamma$ such that $\ol{P_n(x)}\cap J(f_a)$ is connected for all $n$. Moreover,
\begin{enumerate}
\item either the intersection $\bigcap \ol{P_n(x)}$ reduces to $\{x\}$, 
\item or the end of the critical point is periodic : $\exists k>0$  such that for all $n$ large enough $f^k :P_{n+k}(-a)\to P_n(-a)$ is quadratic like. It follows that  the map is renormalizable : there exist $c\in \C$ and  a quasi-conformal homeomorphism  $\phi :P_n(-a)\to V$ , where $V$ is a neighborhood of the filled Julia set of $Q_c$ and $\phi$ conjugates the maps where it is defined. Moreover, for any $x\in J(f_a)$  either the impression  $\bigcap\ol{P_n(x)}$ reduces to $\{x\}$ or to an iterated  pre-image of the critical impression  $\bigcap\ol{P_n(-a)}=I(-a)$.
\end{enumerate}}
\noindent {\sc Proof of the Claim } : This result follows from the construction of the graph done in~\cite{Ro1} that we recall now.
 Denote by $B$ the immediate basin of attraction of $0$. 
The graph $\Gamma$ under consideration is the union of two cycles of rays and two equipotentials. 
More precisely,  we take in $B$ the cycle  generated by the internal  ray of angle $\theta$  of the form $\frac{1}{d^l-1}$ (for any $l$ large enough), the landing point is a repelling cycle, we then  take  the  cycle of external rays landing at this repelling cycle (on the boundary of $B$), for the equipotentials we take any   internal equipotential (in $B$) and any  external equipotential. It is not difficult to see (compare~\cite{Ro1}) that any point of  the Julia set lying in some  sector $U(\theta,\theta')$ (defined below)  is surrounded by a non-degenerate annulus  of the from $P_n(x)\setminus \ol{P_{n+1}(x)}$ ({\it i.e.} lies   in the central component of such an  annulus). The sector $U(\theta,\theta')$ (with $\theta'<\theta$) is defined as follows. Consider the curve  $C$ formed by the internal  rays in $B$ of angles $\theta/d$ and $\theta' +1/d$ and the external rays landing at the corresponding point, $U(\theta, \theta')$  is the connected component of the complement of this curve in $\C$ that contains  the internal ray of angle $0$.
  We need to show that any point of the Julia set will fall under iteration in this domain, then  using Yoccoz's Theorem (see~\cite{Ro1}) we will  get the announced Claim. 
 Any point 
of the filled Julia set belongs to a limb and limbs are sent to limbs (except for the critical limb).
 Therefore,  any point not in the critical limb will fall under iteration in $U(\theta,\theta')$.
 If the critical limb is fixed then it is attached by a fixed ray and already belongs to $U(\theta,\theta')$.
 Otherwise, we look at the limb of the critical value and its orbit will fall in $U(\theta,\theta')$ since the angle of the critical limb is necessarily periodic (indeed the sector of the wake containing the critical value  has angular opening multiplied by $d$ as long as it is not in the wake containing the critical point).

\vskip 0.5em 
 As a direct consequence of Yoccoz's result,  if the map is not renormalizable, the whole Julia set is locally-connected. We will now  consider the case where $f_a$ 
 is renormalizable. 
\vskip 0.5em 
$*$ {\it  If the map $f_a$ is renormalizable then we are in case $2$ of the Claim~:}
Let $K$ denote the filled Julia set of the renormalization $f_a^k$ containing the critical point $-a$. 
We can choose $\theta$ such that the graph previously constructed  does not cut $K$. Indeed,  the graph is forward invariant and any intersection point between $\Gamma$ and  $K$ would  be iterated to a point of  the periodic cycle on the boundary of $B$. Hence, if $K$ intersects $\partial B$ under a cycle,  it is enough to choose $\theta$ of a different period. Now, every  puzzle piece $P_n(-a)$ contains the entire set $K$. Moreover,  since  $K$ is periodic of some period $k$, the puzzle pieces $P_n(-a)$ are all mapped by $f^k$ to $P_{n-k}(-a)$  as a quadratic like map (since $\ol {P_{n+1}(-a)}\subset P_n(-a)$ for large $n$ by the proof of the Claim) so the critical point doesn't escape. Therefore we are in case $2$ of the Claim.
\vskip 0.5em 
$*$ {\it Now we  can assume that  $K=I(-a)$ by   taking the renormalization of lowest period. The  filled in Julia set $K(f_a)$ is the union of $K$ and ``limbs" of it.}

There are two rays landing  at the non-separating fixed point (called $p$) of $f_a^k$ in $K$ (the fixed point corresponding by the conjugacy to the $\beta$ fixed point) and exactly two rays landing at the pre-image of $p$  by $f^k_{\vert K}$ (preimage in $K$). Indeed, if there would be more than two rays landing at  $p$, then they define  some new  sector invariant by $f_a^k$ and it should contain some part of   $K$, which gives the contradiction (the point $p$ would be separating). These two rays  separate $K(f_a)$ in three components, one containing  $B$, denoted by $L$,  one containing neither $B$ nor $K$, denoted by $L'$.   The iterated pre-images of $L'$ by $f_a^k$ and $L$ are called  the limbs of $K$.
 A connected component of $K(f_a)\setminus K$ except $L$ is mapped to a connected component of $K(f_a)\setminus K$. Therefore, any limb different from  $L$ is  an (iterated) pre-image of $L'$. 
\vskip 0.5em 
$*$ {\it For any point $x\in K$, for any neighborhood $U$ of $x$, there exist a sub-neighborhood $V\subset U$ such that $U\cap J(f_a)$ has finitely many connected components. Therefore $J(f_a)$ is locally-connected at the points of $K$. }

Since $K(Q_c)$ is locally-connected, the image $K$ is also locally-connected. Therefore,  there exist a  neighborhood $V\subset U$ such that $V\cap K$ is connected. We  prove now that the diameter of the limbs of $K$ tends to $0$ so that only  finitely many of them 
enter $V$  without being totally included in $V$. 
For this purpose we prove that the diameter of $f_a^{-n}(L')$ tends to $0$, meaning that for any $\epsilon>0$ only finitely many of them have diameter greater than $\epsilon$. 
For this, we shall use Yoccoz's puzzle  for $K(Q_c)$. This puzzle is defined when the non-separating fixed point $\alpha(Q_c)$ is repelling ({\it i.e} when both fixed points are repelling). Therefore we first consider  the case when the fixed point $\alpha$ is not repelling. If it is attracting or parabolic
then the map $f$ is geometrically finite and the result follows from~\cite{TY}. If the point is an irrationally indifferent  fixed point
for $Q_c$, then  its image by  the conjugacy will have the same rotation number (see~\cite{Na}). In the Cremer case both Julia sets (are at the same time of Cremer type and) are non-locally-connected. In  the Siegel case, if the critical point is not on the boundary of the Siegel disk for one map so it is for the other by the conjugacy, and both  Julia sets are non-locally-connected.
 If the critical points are on the boundary of the respective Siegel disks then the post-critical sets  stay in the boundary of the Siegel disks and remain away from the other fixed point and their first pre-images. 
Therefore  $\ol {L'}\cap \ol {\cup_{j\ge 0}f_a^j(-a)}=\emptyset$ and  then using the so called shrinking Lemma (expansion in $\C\setminus \ol {\cup_{j\ge 0}(f_a^j(-a)\cup f_a^j(0) )}$), the diameter  of $f_a^{-n}(L')$  goes to $0$. 

Now we consider the case where the $\alpha$ fixed point of $Q_c$ is repelling. The graph $\Gamma_0$ defining the Yoccoz's puzzle for $Q_c$ is the union of the external rays landing at the point $\alpha$ and some external equipotential. Let us define a new graph $\tilde \Gamma$ for $f_a$ which is a combination of $\Gamma$ and the cycle of external rays landing at the image of $\alpha$ by the conjugacy $\phi$. The puzzle pieces have the same combinatorics for  the map $f_a$ and for the map $Q_c$ using the conjugacy $\phi$ that allows to identify puzzle pieces.  
Two cases appear in Yoccoz's result : either the map $Q_c$ is non-renormalizable  and then the nest of puzzle pieces shrink to points, or it is renormalizable and then it is easy to see that   the orbit of $0$---the critical point of $Q_c$---is bounded away from the $\beta$  fixed point and  its pre-image $-\beta$. In the second case, using the conjugacy $\phi$ one obtains the result  applying the shrinking lemma since the  post-critical set of $f_a$  will be disjoint from the limb $L'$. 
In the first  case we get  a sufficiently small neighborhood of $x$ such that  the intersection  with $J(f_a)$ is connected since the diameter of the puzzle pieces in a nest for $f_a$ shrinks to $0$ also.

\vskip 0.5em 
$*$ {\it Conversely we assume now that the Julia set of $f_a$ is locally-connected. We  prove that the Julia set $K$ of any renormalization of $f_a$ around $-a$ is also locally-connected.}
 
Let $\Phi$ denote the Riemann map of the complement of $K(f_a)$ (which in fact coincides with the B\"ottcher coordinate).
Then $\Psi=\Phi^{-1}$   extends continuously to the boundary. The pre-image $K'=\Psi(K)$ is a compact subset of  the unit circle. 
Therefore its complement is a countable union  of open intervals in the unit circle. Let $\Pi$ be the projection from the unit circle to 
itself that collapses those open intervals to points, {\it i.e.} identifies the whole interval to one point.  
If $t,t'$ are  boundary points of such an open interval, then the external rays of angle $t,t'$ land at the same point in $K$. 
Indeed, these two landing points are in $K$ and the landing point of any external ray in the interval between $t$ and $t'$ is not in $K$ which is a connected set.
  Then we can define a map from the unit circle to $K$ as follows. For $\theta=\Pi(t)$, define $\overline{\Psi}(\theta)=\Psi(t)$.
 By the previous discussion,  this map $\overline{\Psi}$ is well defined and continuous. Therefore $K$ is locally-connected as the continuous image of the unit circle.
\cqfd

 \section{Rational maps}
 {\it  Which Fatou  components of a rational  map are  Jordan domains? }

The property of having a bounded  Fatou component has no meaning for  a rational map. 
The question is  which properties  of the bounded Fatou components of polynomials are used in the proofs of local connectivity results. 
Before we consider  this issue, it is natural to ask if there exist rational maps with  connected Julia sets but with   Fatou components whose  boundaries are not locally-connected. We should consider only rational maps that are sufficiently far from  polynomials. 

 \begin{definition}\label{veritable}
We  say that a rational map is  {\it veritable} if it is not topologically conjugate to a polynomial on its Julia set.
\end{definition}
Notice that in~\cite{Ro2} we introduced the notion of  a  {\it genuine rational map} which by definition is a rational map  that is 
 not conjugate to a polynomial in a neighborhood of its Julia set. This condition is stronger than Definition~\ref{veritable}. Indeed, rational maps of degree~$2$ with a fixed parabolic point at infinity of multiplier $1$  are conjugate to quadratic polynomials on their Julia set (except in some special cases) but cannot be conjugate on a   {\it neighborhood} of their Julia set simply because of the presence of a parabolic basin in the Fatou set (see~\cite{PR}). 

\subsection{Some rational maps as examples.}
$*$ {\it Positive results }(see~\cite{Ro4})
\vskip 0.5em
 Let us start with rational maps of low degree that fix a Fatou component containing a critical point.
 When the rational map is of degree $2$ and the basin is attracting, the map  is necessarily conjugate on its Julia set to a polynomial of degree $2$\,; on can easily see it 
 by using a surgery procedure.  In the parabolic case,  using McMullen's result  (\cite{aut}) the rational map is conjugate to the one that fix infinity with multiplier $1$ and as we mentioned above those maps are conjugate on their Julia set to a quadratic polynomial if the Julia set does not contain 
a fixed Cremer point or a fixed Siegel disk (see~\cite{PR}). Therefore, the question of local connectivity 
for the boundary of the basin  is almost equivalent to the same question for quadratic polynomials.

A rational map of degree~$3$ has $4$ critical points in  $\widehat \C$. 
First assume that three of them are fixed. It is then easy to see  that the rational map is conjugate (by a Moebius transformation)  to a  {\it Newton method } associated to a  polynomial of  degree  $3$, {\it i.e.} to the  rational map   $N_P(z)=z-P(z)/P'(z)$ where $P$ is a cubic polynomial 
with  distinct roots. One should notice that if the fourth critical point is in the immediate basin of attraction of one of the three fixed points,  then
 $N_P$ is conjugate in a neighborhood of its Julia set to a polynomial of degree $3$.
 
In  other cases, the  fixed Fatou components are always Jordan domains (see~\cite{Ro4})\,:
\begin{thm}Let  $N$ be  a cubic Newton method that is a veritable rational map. Then the Fatou components containing a critical point are   Jordan domains. Moreover, the Julia set is locally-connected as soon as 
there is no ``non-renormalizable'' Cremer or Siegel point. 
\end{thm}
  
Let us  now consider the set of rational maps of degree $2$ having a cycle of period $2$ of Fatou components containing a critical point. One example is the family of rational maps with a period $2$ critical point studied in~\cite{AY}. In this article the authors prove that if   the map is  non-renormalizable, with only repelling periodic points, then it has a locally-connected Julia set. 
It  seems  reasonable to believe that this result holds in general in this family and that it persists when the cycle becomes  parabolic.

One may even wonder  whether  the critical  Fatou components of a veritable rational  map would always have  locally-connected boundaries.
\vskip 1.5em
\noindent $*$  {\it Negative results }(see~\cite{Ro2})
\vskip 0.5 em
In~\cite{Ro2} one exhibits two families of examples  that illustrate the following result\,:
\begin{thm}\label{th:nonlc} There exist veritable rational maps with connected Julia sets that posses a Fatou component with non-locally-connected boundary.
\end{thm}
The first set of examples can be found in the works of   Ghys and  Herman~\cite{Ghys,Herman}\,; they are in  the family
  $\displaystyle f_{a,t}=e^{2i\pi t}z^2\frac{z-a}{1-az}$ ($a>3$).
  
 The second set  of examples can be found  in the family 
   $\displaystyle g_a(z)=z^3\frac{z-a}{1-az}$, $a\in \C$. They are   obtained from perturbing  a map $g_{a_0}$ with a parabolic point in order to create a Cremer point.  

We will briefly explain in next section why these maps have a  non-locally-connected Julia set. Nevertheless, such examples are not satisfactory since, roughly speaking,  one can  see the trace of a Julia set of a quadratic polynomial in them.

\section{Douady-Sullivan criterion}\label{s:conj}
All our examples of rational maps and polynomials with a Fatou component whose boundary is not locally-connected, share a certain property.
 We will call it the {\it   Douady-Sullivan criterion} since it was originally used   by them 
to prove that the Julia set of a polynomial with a Cremer point is not locally-connected.

First note that if  a periodic Fatou component  $B$ is not  simply connected, then  $\partial B$ is not connected and $J(f)$ is not connected
 (one can also deduce  that  $\partial B$ is not locally-connected at any point of its boundary). 

\begin{definition}[(Douady-Sullivan criterion)]\label{critere de douady-sullivan} A rational map 
 $f$ is said to satisfy the   {\it   Douady-Sullivan criterion}
  whenever 
 $f$ has a $k$-periodic   Fatou component  $B$ that is simply connected and  contains a critical point, and if there exists 
   compact  set $C $  in the boundary of $B$ such that
  \begin{itemize}
  \item   $C$ does not contain any critical point of  $f^k$ and 
  \item   the  restriction  $f^k:C\to C$ is a bijection.
\end{itemize}
\end{definition}
\begin{lem}\label{l:DS} Let $f$  be a  rational map that satisfies the {\it   Douady-Sullivan criterion}.
Let $C$ denote the compact set and $B(p)$ the Fatou component appearing in the definition of the criterion, 
 with $C\subset \partial B(p)$. Then \;:
\begin{itemize}
\item
either $\partial B(p)$ is not  locally-connected, 
  \item or 
  $C$
is the finite union of parabolic or repelling cycles.
\end{itemize}\end{lem}
\proof (Compare~\cite{Mi1}). 
The basin  $B(p)$ is simply connected, assume that  its  boundary is   locally-connected. We assume that $k=1$  replacing    $f^k$ by $f$.  From Carath\'eodory's Theorem, we know that the map is conjugate on the boundary to one of the models  $Z_d(z)=z^d$ or $B_d(z)=\frac{z^d+v_d}{1+v_dz^d}$ (where $v_d=\frac{d-1}{d+1}$) on the closed unit disk. Notice that the restriction of map $B_d$ to the unit circle is  topologically conjugate to the restriction of  $Z_d$ to the unit circle. Therefore there exists a map    $\gamma:\S^1\to\partial B(p)$ 
that is a  semi-conjugacy between     $e^{2i\pi\theta}\  \mapsto e^{2i\pi d \theta}$ and $f$.
One considers then the set  $\Theta=\{\theta\in \R/\Z\mid \gamma(e^{2i\pi\theta})\in C\}$.
The map $m_d:\theta\mapsto d\theta$ is a bijection from  $\Theta$ to itself.
Indeed, if two rays landing  at the same point are mapped onto the same ray, then they land at a critical point, which is excluded by our hypothesis in Definition~\ref{critere de douady-sullivan}. Moreover, every point of $C$ has a pre-image in $C$ and  at least  one ray  lands at this point\;;  therefore $m_d:\Theta\to \Theta$ is surjective. It follows then that $m_d:\Theta\to \Theta$ is a 
 homeomorphism since $\Theta$ is compact. Finally notice that, $m_d$
 is expanding and hence  $\Theta$ must be   finite. 
Indeed,  cover  $\Theta$ by a finite number  $N$ of balls of radius   $\epsilon$ sufficiently small\;;  since  $m_d$ is a homeomorphism, the  pre-image of a ball of radius $\epsilon$ is a ball  of radius  $\epsilon/d$, so  $\Theta$ is covered by 
 the union of  $N$ balls of diameter  $\epsilon/d$, etc. It follows that,   $\Theta$  is the  union of $N$ points and those points  are pre-periodic angles. Since  $m_d$ is a  bijection,   $\Theta$  is a union of periodic 
  cycles.
    One deduces that  $C$ is the union  of   cycles of $f$. These cycles are   parabolic or repelling by the Snail Lemma (see~\cite{Mi1}).
\cqfd

\subsection{Douady-Sullivan criterion  in the previous examples.}\label{s:CremerSiegel}
\noindent
{\bf Polynomial with a  Cremer point.}
 
Let  $f$ be a  rational map. Recall that a {\it Cremer point } of  
 $f$ is a point of the Julia set   $J(f)$ that is  irrationnally indifferent.  The  Julia  set of a polynomial  with a   Cremer point is not locally-
connected. This follows from  Lemma~\ref{l:DS}   since  $f$ satisfies the 
  Douady-Sullivan criterion with $p=\infty$ and   $C$ being the cycle generated by the Cremer point.

 \vskip 0.5em
 \noindent{\bf Polynomials with  Siegel disks.}
 
 Let  $f$ 
be a polynomial and $\Delta$ a periodic Siegel  disc for $f$
 such that no critical point is on the  cycle generated by the boundary $\partial \Delta$. Then the Julia set  $J(f)$ is not 
locally-connected  (see~\cite{Mi1} for instance).
Indeed,     $f$ satisfies the 
  Douady-Sullivan  criterion taking $p=\infty$ and   $C$ to be the cycle generated by boundary   $\partial \Delta$.

 \vskip 0.5em
 \noindent{\bf Rational maps\,: example of Ghys-Herman.}
\vskip 0.5em
Now we consider  the first  family of examples studied in~\cite{Ro2}, namely the family $\displaystyle f_{a,t}=e^{2i\pi t}z^2\frac{z-a}{1-az}$ with  $a>3$ and $t\in \R$. The restriction of $f_{a,t}$  to  $\S^1$ is an  $\R$-analytic diffeomorphism. According to  Denjoy's Theorem, if the rotation number  $\alpha= \rho(f_{a,t})$ is irrational,  $f_{a,t}$ is topologically  conjugate on $\S^1$ to the    rigid rotation $R_\alpha$ by some  homeomorphism $h_{\alpha,t}$. 
  E.~Ghys shows  (in~\cite{Ghys}) that if   $h_{\alpha,t}$ is  quasi-symmetric but not   $\R$-analytic, then the 
  polynomial $P_\alpha(z)=e^{2i\pi\alpha}z+z^2$ has a Siegel disk whose  boundary is a quasi-circle  not containing the critical point. 
  On the other hand, to compare   $f_{a,t}$ with  $P_\alpha$ Ghys performs a surgery that provides a
 homeomorphism $\psi$ such that $\psi(\S^1)$ is the boundary of the  Siegel disk and such that the boundary of the immediate basin of  $\infty$ for  $f_{a,t}$ is the  image by  $\psi$ of  $J(P_\alpha)$. 
 As noticed previously, the Douady-Sullivan  criterion  implies that the  Julia set of  $P_\alpha$ is not 
 locally-connected since the boundary of the Siegel disk contains no critical point.
 This implies that the boundary of the basin of $\infty$ for $f_{a,t}$ is not locally-connected.
 To conclude, we use the following result of    M.~Herman\;:  for any   
$ a>3$, there  exists values of $ t\in \R$    such that the conjugacy $h_{a,t}$ between
 $f_{a,t}$ and  $R_\alpha$ is  quasi-symetric but not  $C^2$. \vskip 0.5 em 
 
\begin{figure}  \begin{center}\vskip -10em 
\includegraphics[height=6 in]{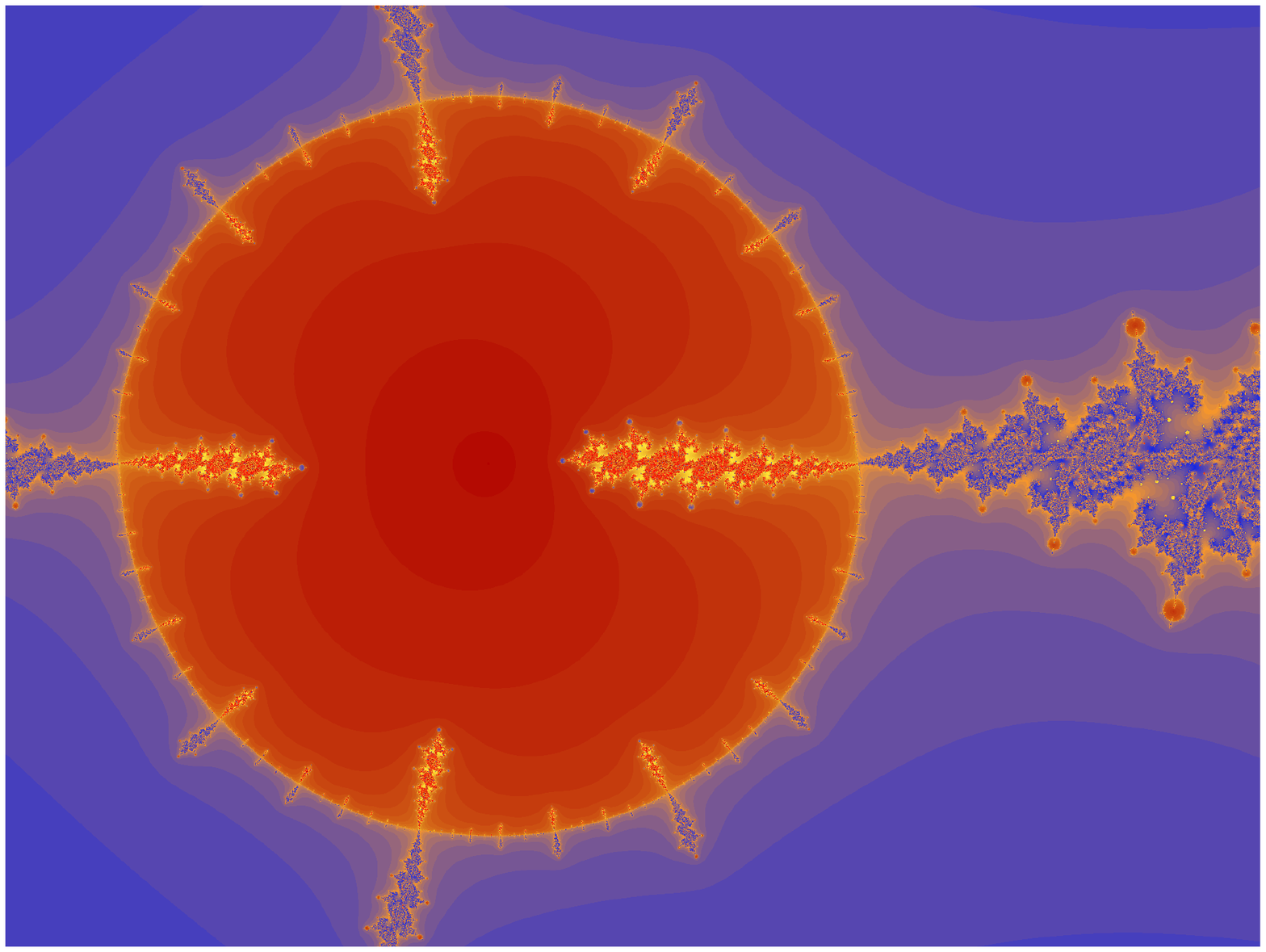}
\vskip -10em \caption{ \it The Julia set of  $g_a$ for a value  $a$ next to  $5$.}
\end{center}
\end{figure}

 \begin{figure}  \begin{center}\vskip -10em 
\includegraphics[height=5 in]{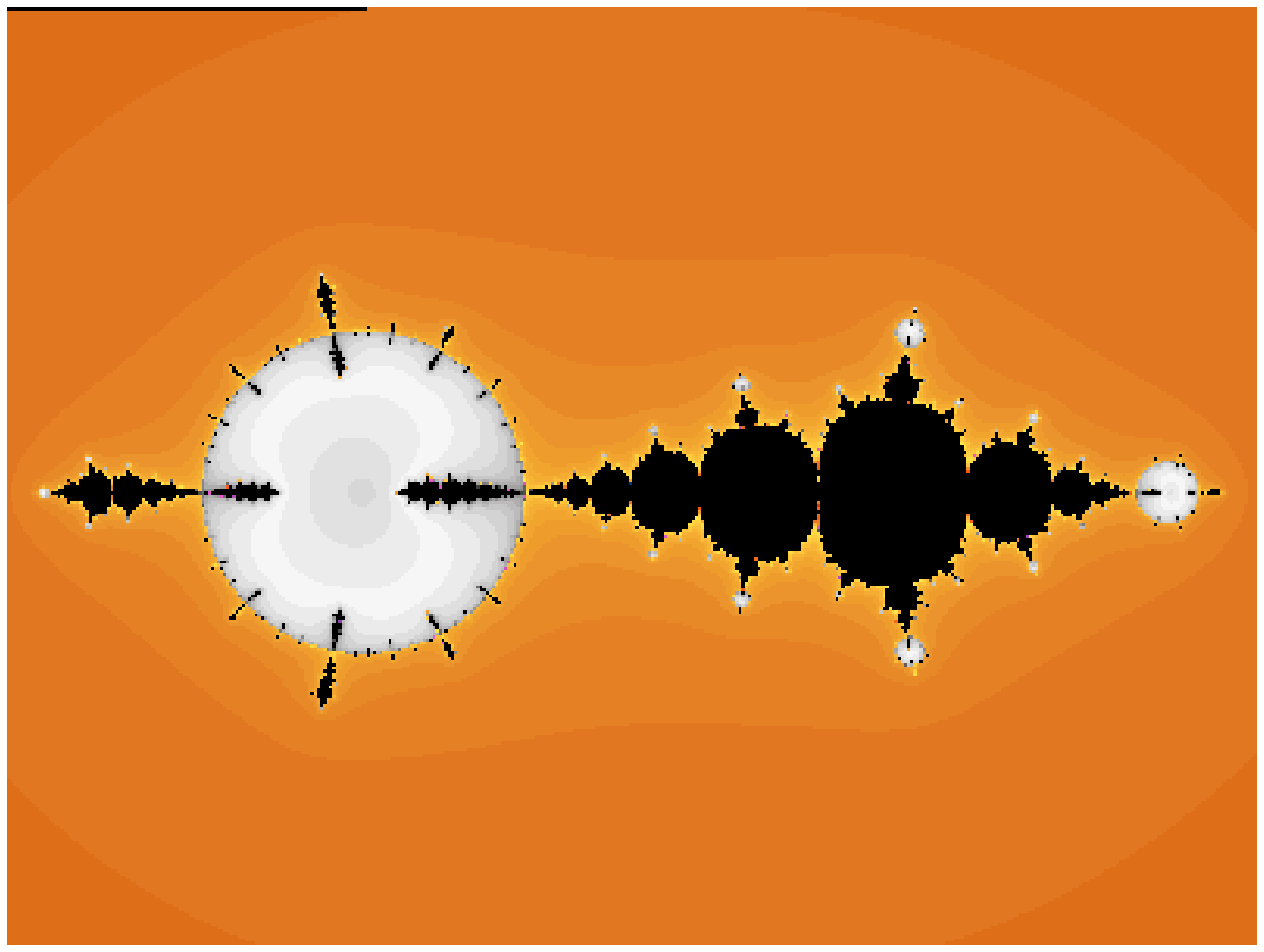} 
\vskip -8em \caption{ \it The Julia set of  $g_5$.}
\end{center}
\end{figure}\vskip 0.5em
 \noindent{\bf Rational maps\,: perturbation of a  fraction of  Blaschke.}
 
 In the second family of  examples studied in~\cite{Ro2}, namely   { $\displaystyle g_a(z)=z^3\frac{z-a}{1-az}$} with  $a\in\C$,  the proof is much easier.
 One sees directly that the  Blaschke product $g_5$  is  renormalizable. Indeed, it admits a restriction that is quadratic like in some open set bounded by rays in the immediate basin of the attracting fixed point $0$. This  restriction admits  a parabolic   point at  $2-\sqrt 3$.
There exists a neighborhood of  $a_0=5$ in the parameter space such that  for $a$ in this neighborhood the map $f_a$  has  a fixed point  $p(a)$ that is a   holomorphic function satisfying 
  $p(5)=2 -\sqrt 3$.  Moreover in this neighborhood,  $g_a$ admits a restriction that is  polynomial-like in the neighborhood of  $p(a)$.
  Therefore one can find values of $a$ near  $a_0$ such that 
  $g_a$ is renormalizable,  with the renormalized  filled Julia set containing a Cremer point. 
  Since the open sets defining the renormalization intersect the immediate basin of $0$, the Cremer point thus obtained has to be on the boundary of the immediate basin. 
  Those maps  $g_a$ verify the   Douday-Sullivan  criterion just by taking for  compact set $C$ the Cremer point.

\subsection{Infinitely renormalizable polynomials}\label{s:infinitelyDS}

Finally there is a class of examples which we have not yet discussed.
Indeed, in the class of infinitely  renormalizable polynomials,
 one can find polynomials having connected but not locally-connected Julia set.
Several works have been devoted to their studies, see for instance  \cite{So,M3,MiNlc,Levin-2009}.

We consider in this section particular infinitely renormalizable polynomials\,; these are the quadratic polynomials $Q_c(z)=z^2+c$ where $c$ is in a  \textit{limit of a sequence $(H_n)$ of hyperbolic components such that $H_{n+1}$ is attached to $H_n$ with higher period}. Let us be more precise.

\begin{definition}
\begin{itemize}
\item We will use the notation $\mandelbrot$ for the classical Mandelbrot set, that is, the set of $c\in\complexes$
 such that the orbit of the critical point of $Q_c$ is bounded.

\item Let $H$ be a hyperbolic component of $\mandelbrot$, {\it i.e.}, a connected component of the interior of $\mandelbrot$ such that $Q_c$  has an attracting periodic point of some period $k$ for every $c\in H$.
In $H$, there exists a unique parameter $c\in H$ such that the critical point of the quadratic polynomial $Q_c$ is periodic. We call the period of this point  \textit{the period of} $H$ and this parameter $c$ is called \textit{the center of} $H$.

\item We say that a hyperbolic component $H'$ of $\mandelbrot$ is \textit{attached to}
 $H$ if its boundary intersects the boundary of $H$.
In this case, their boundaries intersect at a unique  point.

This intersection point is called the \textit{root of} $H'$ 
 if the period of $H'$  is greater than the period of $H$.
Every hyperbolic component of $\mandelbrot$
 has at most one root\footnote{And they all have at least one root.
 But this root does not necessarily belong to the boundary of another hyperbolic component.}.
We will write $r(H)$ for the root of the hyperbolic component $H$.
 At the parameter $c=r(H)$, $Q_c$ has a parabolic cycle of period  less or equal to the period of $H$.
 If $Q_c$ has a parabolic point of period equal to the period of $H$, with $c\in \partial H$, we call this the root of the hyperbolic component.
 As we will see later it is the unique point where the multiplier function $\mu:\ol H\to \ol\D$ takes the value $1$. 

\item Let $(H_n)_{n\geq 0}$ be a sequence of hyperbolic components of $\mandelbrot$.
 We say that $(H_n)_n$ is a \textit{chain of  components arising from} $H_0$
 if for all $n\geq 0$, $H_{n+1}$ is attached to $H_n$ and the period of $H_{n+1}$ is greater than the period of $H_n$.
\end{itemize}
\end{definition}

One has to  notice that when $H'$ is attached to $H$ at its root point $r(H')$, the period of $H$ divides the period of $H'$.

The parameters $c$ we will  consider here  are limits  of the sequences of $(r(H_n))_n$ where $(H_n)_n$
 is a chain of components.

When a parameter $c$ belongs to a hyperbolic component $H$ of period $k$,
 the mapping $Q_c$ has an attracting cycle $Z_H(c)$ of period $k$.
 The points in this cycle are holomorphic functions of  $c\in H$.
These holomorphic mappings can be extended to some regions containing the component $H$ but they do not extend to any neighborhood of the root $r(H)$.

 Let $H'$ be a hyperbolic component which is attached to $H$ and whose period is $k'>k$.
The cycle $Z_H(c)$ of $Q_c$, attracting when $c\in H$, becomes  parabolic when $c=r(H')$.
 Then, when $c$ enters $H'$ a bifurcation occurs\,: The cycle $Z_H(c)$ becomes repelling while an attracting cycle
 $Z_{H'}(c)$ of period $k'$ appears.

Let $(H_n)_n$ be a chain of  components arising from $H_0$. Let $k_n$ denote the period of $H_n$.
 It can be proved that, for all $n\geq 0$,
 the  mappings $Z_{H_n}$ have well defined analytic continuations to some neighborhood of 
 $\ds{\bigcup_{m\geq n}}\left.\barre{H_m}\,\right\backslash\{r(H_{n})\}$.

 We denote these continuations by the same $Z_{H_n}$. 
 Note that 
the cycle $Z_{H_n}(c)$ is repelling for
 $c\in\ds{\bigcup_{m\geq n+1}}\left.\barre{H_m}\,\right\backslash\{r(H_{n+1})\}$.

\begin{lem}\label{chaine}
There exist chains of components $(H_n)_{n\geq 0}$ such that
the sequence $(r(H_n))_n$ converges and such that the limit $c_*$ has the following properties :
\begin{enumerate}
\item
for all $n\geq 0$,
the cycle $Z_{H_n}(c)$ converges to a repelling cycle $Z_{H_n}(c_*)$ as $c\tend c_*$\;;
\item
the closed set $Z(c_*)=\ol{\cup_{n\ge 1} Z_{H_n}(c_*)}$  does not contain $0$\;;
\item  $Q_{c_*}$ is infinitely renormalizable and
\item the Julia set $J(c_*)$ of $Q_{c_*}$ is not locally-connected.
\end{enumerate}
\end{lem}

The difficulty in the choice of the sequence $(H_n)_n$
 is to ensure that the distance between the critical point $0$ and the cycles $Z_{H_n}(c)$, for
$c\in \ds{\bigcup _{m\geq n+1}\ol{H_m}}$,
is bounded  below uniformly in $n$.

In \cite{So,M3,MiNlc,Levin-2009} one can find quantitative conditions in terms of the roots $r(H_n)$ which ensure that a chain of  components $(H_n)_n$
 converges to a unique parameter $c_*$ with the above properties (we will come back to this in section~\ref{s:infinitely}).
But if one simply wants to show the existence of a sequence $(H_n)$ satisfying the conditions of  Lemma \ref{chaine}, one can proceed the following way.

We begin by choosing a component $H_0$ of period $k_0$ and a component $H_1$ of period $k_1=rk_0$, attached to $H_0$.
When the parameter $c$ is the root of $H_1$, the cycle  $Z_{H_0}(c)$ is parabolic (with multiplier $\neq 1$)
and it can be followed in a neighbourhood of $r(H_1)$ in $\C$ as a cycle of period $k_0$.
By continuity there exists $\eps >0$ and a neighbourhood $V_1$
of $r(H_1)$ such that for $c\in V_1$ the cycle stays at a distance $\geq \eps$ from the critical point $0$ of $Q_c$.
When $c$ enters $H_1$, the cycle $Z_{H_0}(c)$ becomes repelling and the attracting $k_1$-periodic cycle $Z_{H_1}(c)$ arises.
This cycle appears in $k_0$ clusters of $k_1/k_0$ points (called \textit{bifurcated cycles}) around the points of the cycle $Z_{H_0}(c)$.

This description is valid at all points of ${H_1}$ and remains valid on the boundary of the component $H_1$.
 It is easy to see that the distance between the cycle $Z_{H_1}(c)$ and the cycle $Z_{H_0}(c)$ tends to $0$ as $c\in \ol{H_1}$ tends to $r(H_1)$.
Let us choose a component $H_2$ which is attached to $H_1$, such that its root
 $r(H_2)$ belongs to $V_1$, and such that the distance between the cycle $Z_{H_1}(r(H_2))$ and $Z_{H_0}(r(H_2))$ is  $< \eps /3$.
By continuity, there exists a neighbourhood $V_2$ of $r(H_2)$ contained in $V_1$ such that for $c\in V_2$, the distance between the cycles
 $Z_{H_1}(c)$ and $Z_{H_0}(c)$ is less than $\eps/3$. Now we choose $V_2$ disjoint from $\ol{H_0}$, so that for all $c\in V_2$ the cycle $Z_{H_0}(c)$ is repelling.
Repeating this argument one can build by induction a chain of components $(H_n)_n$ and a decreasing sequence of neighbourhoods $V_n$ of $r(H_n)$,
 which are disjoint from $\ol{H_{i}}$ for all $i\le n-2$ and such that for all $c\in V_n$ the distance between the cycles
 $Z_{H_{n-1}}(c)$ and $Z_{H_{n-2}}(c)$ is less than ${\eps}/{3^n}$. 

We can also choose the $V_n$ such that  their diameters tend to $0$.
Hence the sequence of the roots $r(H_n)$ converges to a point $c_*$.
It is then easy to check that for all $n\geq 0$, the cycle 
 $Z_{H_n}(c_*)$ is repelling and is at a distance of at least $\eps/2$ from the critical point $0$ of $Q_{c_*}$.

Now we give a proof of the fact that the Julia set of $Q_{c_*}$
is not locally-connected, only using the Douady-Sullivan criterion.

\begin{lem}\label{lem ref pour zq}
 Let $c_*$ be a limit point of a chain of 
 components $(H_n)_n$ such that for all
 $n\geq 0$, the cycle $Z_{H_n}(c)$ converges to a repelling cycle $Z_{H_n}(c_*)$
 and such that the distance between the critical point $0$
 and the cycles $Z_{H_n}(c)$, when $c\in \bigcup _{m\geq n+1}{H_m}$,  
 is bounded below
 uniformly in $n$. Then $Q_{c_*}$ satisfies the Douady-Sullivan criterion.
\end{lem}

\proof

For $q\geq 2$, we denote by $\relatifs_{q}$ the group of integers modulo $q$.
In order to lighten notations,
 if $i\in\{0,\dots,q-1\}$, we will use the number $i$ for its residue class modulo $q$.

We consider the compact set $C=\ol{\bigcup _{n\geq 0}Z_{H_n}(c_*)}$, which by the  assumption  does not contain the critical point.
 Let us  show that the mapping $Q_{c_*}:C\to C$ is bijective.

The surjectivity of $Q_{c_*}$  on $C$ follows from the fact that $Q_{c_*}$ is onto on $\bigcup _{n\geq 0}Z_{H_n}(c_*)$ since each $Z_{H_n}(c_*)$ is a cycle.

The injectivity is more subtle, although it is obvious that $Q_{c_*}$
 is injective on $\bigcup _{n\geq 0}Z_{H_n}(c_*)$.
We label the points of the cycles along clusters.
Let $k_i$ be the period of $H_i$. Define the numbers $q_i$ for $i\geq 0$ by the relation $k_{i+1}=q_ik_{i}$.
For  later use we set $q_{-1}=k_0$.

Let $z_{i}$ denote, for $ i\in \relatifs_{k_0}$, the points of the first cycle $Z_{H_0}(c)$
 in such a way that $Q_c(z_i)=z_{i+1}$ (where the indices are taken  modulo $k_0$).
Then, for $i\in\relatifs_{k_0}$ and $j\in\relatifs_{q_0}$, we let $z_{i,j}$ be the part of the cycle $Z_{H_1}(c)$
which bifurcates from the point $z_i$, with indices chosen so that,
 if $ i\neq k_0-1$ then $Q_c(z_{i,j})=z_{i+1,j}$ and if $i=k_{0}-1$ then $Q_c(z_{i,j})=z_{0,j+1}$.

For any $n\geq 0$, the cycle $Z_{H_{n+1}}(c)$ bifurcates from the cycle $Z_{H_{n}}(c)$.
Because of this, we can label the points $z_{\eps_0,\cdots,\eps_n}$ of the cycle $Z_{H_{n}}(c)$
 according to the dynamics.
More precisely, there is a mapping $\tau_n$ from
$\relatifs_{k_0}\times\relatifs_{q_0}\times\cdots \times\relatifs_{q_{n-1}}
=\relatifs_{q_{-1}}\times\relatifs_{q_0}\times\cdots \times\relatifs_{q_{n-1}}$
 into itself,
such that $Q_c(z_{\eps_0,\cdots,\eps_n})=z_{\tau(\eps_0,\cdots,\eps_n)}$.
The mapping $\tau_n$  is defined in the following way.
For  $( \eps_0,\cdots,\eps_n)\in\relatifs_{q_{-1}}\times\relatifs_{q_0}\times\cdots \times\relatifs_{q_{n-1}}$, let $(\eps'_0,\cdots,\eps'_n)=\tau( \eps_0,\cdots,\eps_n)$ and 
let $0\leq j\leq n$ be the smallest integer such that $\eps_j\neq k_{j-1}-1$.
 Then, for $l<j$, $\eps'_l=0$, $\eps'_j=\eps_j+1$ and for $i>j$, $\eps'_i=\eps_i$.
Moreover the image of $(q_{-1},\dots,q_{n-1})$ by $\tau$ is $(0,\dots,0)$.

Taking the limit as $n\tend \infty$, this definition yields a mapping $\tau$
 from $\ds{\prod_{i=0}^\infty} \relatifs_{q_{n-1}}$ into itself which is bijective.
Now we need to check that the  ``parametrization''
  $\ueps
  \in\ds{\prod_{i=-1}^{+\infty}}\relatifs_{q_{i}}
   \mapsto z_\ueps\in C$ is injective.
Without loss of generality one can assume that the sequence of hyperbolic components arises from the main cardioid of $\mandelbrot$.
In order to find neighborhoods which group clusters of bifurcated cycles together
 we define disjoint graphs $\Gamma_1(c), \dots, \Gamma_n(c),\dots$ satisfying the following properties (see~\cite{Rhabi}):
 \begin{itemize}
\item $\Gamma_n(r(H_n))$ is made of external rays landing at the parabolic cycle of $Q_{r(H_n)}$\,;
\item $\Gamma_n(c)$ exists in a neighborhood $U_n$ of $\ds{\ol {\bigcup_{i\ge n}  H_i}}\;\sans\{r(H_n)\}$\,;
\item $\Gamma_n(c)$ depends continuously on $c$ inside $U_n\cup \{r(H_n)\}$\,;
\item $\ds{\bigcup_{k=1}^n}\Gamma_k(c)$ separates the points of the cycles  $Z_{H_n}(c)$ but not the points of the cycle 
  $Z_{H_m}(c)$ where $m>n$\,;
\item  $\ds{\bigcup_{k=1}^n}\Gamma_k(c)$  separates the points $z_{\eps_0,\cdots,\eps_n\cdots}$
 which differ in at least one term $\eps_i$ for $0\le i\le n$.
\end{itemize}

We obtain the graph $\Gamma_n(r(H_n))$
 by considering the cycle of external rays landing at the parabolic cycle for the parameter $ c=r(H_n)$.
There exists a holomorphic motion of this graph defined in a region containing $H_{n}$
 and bounded by external rays in $\C\setminus \M$ landing at $r(H_n)$.
At the parameter $r(H_n)$, the rays which $\Gamma_n(r(H_n))$ is made of
 separate the critical points of the iterate $Q_{r(H_n)}^{ k_n}$. It follows that they also separate the point of the attracting cycle $Z_{H_n}$.
Since the rays and the cycles can not cross each other (the period being different),
the graph $\bigcup_{k=1}^n\Gamma_k(c)$
  separates the points $z_{\eps_0,\dots,\eps_n\dots}$ which differ in at least one term $\eps_i$
 for  $0\le i\le n$.
Hence, if $\ueps\neq \ueps'$ then $z_{\ueps} \neq z_{\ueps'}$.
The injectivity for finite sequences is obvious.
 
 Assume  now that two distinct sequences $z_{\epsilon_0,\cdots,\epsilon_n}$ and $z_{\epsilon'_0,\cdots,\epsilon'_n}$ converge respectively to $z$ and $z'$ such that $z\neq z'$.
Then $z=z_{\ueps}$ and $z'=z_{\ueps'}$ with $\ueps \neq \ueps'$.
From $Q_{c_*}(z_{\ueps})=z_{\tau(\ueps)}$ and $Q_{c_*}(z_{\ueps'})=z_{\tau(\ueps')}$ and from the fact that $\tau$ is injective,
 it follows that $Q_{c_*}(z)$ and $Q_{c_*}(z')$ are distinct.
As a consequence, $Q_{c_*}$ is injective on $C$. \cqfd

We can consider the same  question with a sequence of primitive renormalizations. We say that a parameter is {\it primitive renormalizable} if it belongs to a primitive copy of $\M$ in $\M$, {\it i.e.}, contained in a maximal copy of $\M$ in $\M$ which is not attached to the main cardioid of $\M$. We say that a parameter is {\it infinitely primitive renormalizable} if  it belongs to an infinite sequence  $\M_n$, each $\M_n$ being  a primitive copy of $\M$  in $\M_{n-1}$.

\begin{question} Does there exist   infinitely primitive renormalizable quadratic polynomials  having a connected but not locally-connected Julia set? 
\end{question}

\subsection {Conjectures for rational maps.}

With the above examples in mind, the following  conjecture seems   reasonable. 
\begin{conjecture} 
Let $f$ be a rational map whose Julia set is connected. 
If  $f$ has a  periodic  Fatou  component  which contains a  critical point whose boundary is not locally-connected, then   $f$ satisfies the  Douady-Sullivan criterion. 
\end{conjecture}
Notice that in this conjecture rational maps includes polynomials. 

Let us return to the omnipresence of polynomials in our  examples.
 We notice that in each of our examples, the boundary of the periodic critical Fatou component contains a copy of a non-locally-connected quadratic Julia set.
 In the example of Ghys-Herman, the boundary of the   immediate  basin of $\infty$ is   homeomorphic to the Julia set of a quadratic polynomial which is not locally-connected.
 In the degree $4$ Blaschke product example, the boundary of the immediate basin  of $0$ contains  the image
 (by the straightening map of  Douady-Hubbard) of a quadratic Julia set  which is not  locally-connected.
 Motivated by these examples, we propose the following\,:

\begin{conjecture} Let  $f$ be a rational map  whose  Julia set is  connected. 
If  $f$ has a  periodic  critical Fatou  component $U$ whose boundary is not locally-connected, then 
 $\partial U$ contains the homeomorphic image of some non-locally-connected polynomial   Julia set.
\end{conjecture}

Notice that  to be at the boundary of a Fatou component is crucial.
 Indeed, there exist cubic Newton maps $N$ such that the Julia set $J(N)$ contains a quasi-conformal copy of a non-locally-connected quadratic Julia set even though $J(N)$ itself is locally-connected
 (see~\cite{Ro4}).

\begin{figure} [ht] \begin{center}  
\includegraphics[height=4 in]{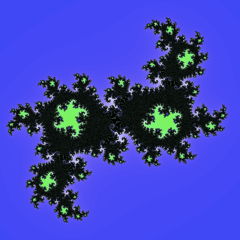} 
\caption{\it The Julia set of some $f_t$. One can imagine the  hairy circle (image courtesy of  H. Inou).}
\label{f:hairs}
\end{center}
\end{figure}

  Notice also that we do not ask that the homeomorphism conjugate the dynamics. 
  Let us consider the map 
$f_t(z)=e^{2i\pi t}z^2(z-4)/(1-4z)$. This is an example in the class of Ghys-Herman studied in section~\ref{s:CremerSiegel}. This map preserves the unit circle, it is of degree $1$, the critical points are not on the unit  circle. Therefore  one can define a rotation number $\rho(f_t)$ of the restriction of the map on the circle.  Since $\rho(f_t)$  is continuous in $t$, one can find some $t$ such that $\rho(f_t)$ is not a Brjuno number. 
This implies in particular that there is no Herman ring around the unit circle. By the theory of Perez-Marco there is  a  ``hedgehog''  with hairs around the circle
 (see figure~\ref{f:hairs}). 
The boundary of the basin of $\infty$ (or $0$) contains this ``hairy circle'' and its pre-image touching at the critical point and all the iterated backward pre-images connected to this.   Can we say that this compact connected set  is homeomorphic to the Julia set of a quadratic polynomial? 
If  that were the case,  the quadratic polynomial   would not be conjugate on its Julia set  to  $f_t$; indeed,   the quadratic polynomial would  have a Siegel disk, but  here the rotation number is non-Brujno.

\section{The case of infinitely satellite renormalizable quadratic polynomials : a model}\label{s:infinitely}
Our aim in this section  is to propose a conjectural condition on some combinatorial data related to an infinitely satellite renormalizable
 quadratic polynomial that implies non-local connectedness its Julia set. The combinatorial data we are interested in
 is the sequence of rotation numbers $(p_n/q_n)_n$ defined in section \ref{sec nombres de rotation}.

We are also interested in a description of the post-critical closure $\ol{P_f}$.
 There are some similarities with the hedgehogs and the Cantor bouquets.

This section contains a description and the beginning of the investigation of a geometric model of the sequence of straightenings of an infinitely satellite
 renormalizable quadratic polynomial which provides such a  conjectural condition.
 We owe the idea of this model to Xavier Buff.

In what follows, $\naturels$ will represents the set of non-negative integers.

\subsection{Combinatorial data for satellite renormalizable polynomials}\label{sec nombres de rotation}

\begin{definition}
 A hyperbolic component $H$ of $\mandelbrot$ is called \textit{satellite to} another hyperbolic component $L$
 if it is attached to $L$ at its root point $r(H)$.
\end{definition}

If $H$ is a hyperbolic component of $\mandelbrot$, the multiplier map $\lambda_H$  of $H$ will refer to
 the mapping that sends a parameter $c\in H$ to  the complex number $\lambda_H(c)$ in $\disque$ that is the multiplier of the unique attracting cycle of $Q_c$.

It is well known that the multiplier map can be extended to a homeomorphism from $\barre{H}$
 onto the closed unit disk.
 Moreover, attached components meet only at parameters at which the multiplier map is a root of unity.

\begin{definition}
 Let $H$ and $L$ be hyperbolic components of $\mandelbrot$ such that $H$ is satellite to $L$.
The rational number $p/q$ such that the multiplier map of $L$ sends the root of $H$ to $e^{2\pi\ii p/q}$
 is called \textit{the rotation number of} $H$ with respect to $L$\,; {\it i.e.} $\lambda_L(r(H))=e^{2\pi\ii p/q}$.
\end{definition}

\begin{definition}
For $c\in\complexes$, the   quadratic polynomial
 $Q_c$ is {\it infinitely satellite renormalizable }  if it is  a limit point of a sequence of 
 hyperbolic components $H_n$ of $\mandelbrot$
 such that $H_{n+1}$ is satellite to $H_n$ for all $ n$. 
   
To each such $c$ and $H_0$ we can associate a 
 sequence of \textit{rotation numbers} $(p_n/q_n)_n$ : it  is the sequence  $p_n/q_n$ of   rotation numbers of $H_{n+1}$ with respect to $H_n$.
\end{definition}

Thanks to the connectedness principle (Theorem~\ref{cp}, compare \cite{mmullen-renorm}),
 if the Julia set of a renormalization  of a polynomial is not locally-connected
 then the Julia set of the original polynomial is not locally-connected.

As a consequence, we are only interested in the tail of the sequence $(p_n/q_n)$
 which is independant of the choice of $H_0$ in the sense that if $\left(H_0,(p_n/q_n)_n\right)$
 and $\left(H'_0,(p'_n/q'_n)_n\right)$
 are both associated to the same parameter $c$,
 then either there is $n$ such that $H'_0=H_n$ and $p'_k/q'_k=p_{n+k}/q_{n+k}$ (for all $k$) or
 $H_0=H_n'$ and $p_k/q_k=p'_{n+k}/q'_{n+k}$.

Keeping this in mind, we will not mention $H_0$
 when we talk about the sequence of rotation numbers of an infinitely satellite renormalizable
 quadratic polynomial.

\subsection{Definition of the model}

Let $(p_n/q_n)_n$ be a sequence of reduced fractions in the interval $]0,1[$, where
 $q_n>0$. We suppose that the sequence $(p_n/q_n)_n$ converges to $0$.

Let $C>1$ be a fixed constant and define $t_n$ as
$$t_n=C\fr{p_n}{q_n}.$$
We refer to the  Lemma \ref{lem inclusions C decroissante}, and the observation following the statement of lemma \ref{lem restriction sur C},  about the role of the constant $C$.

We denote by $M_n$ the M\"obius transformation
$$M_n(z)=\fr{1-t_n/z}{1-t_n}.$$
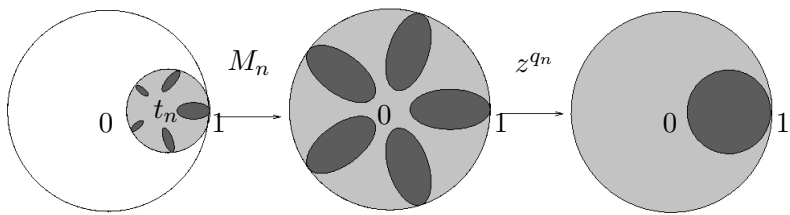
\begin{figure}[ht]
\begin{center}
\scalebox{0.42}{\includegraphics{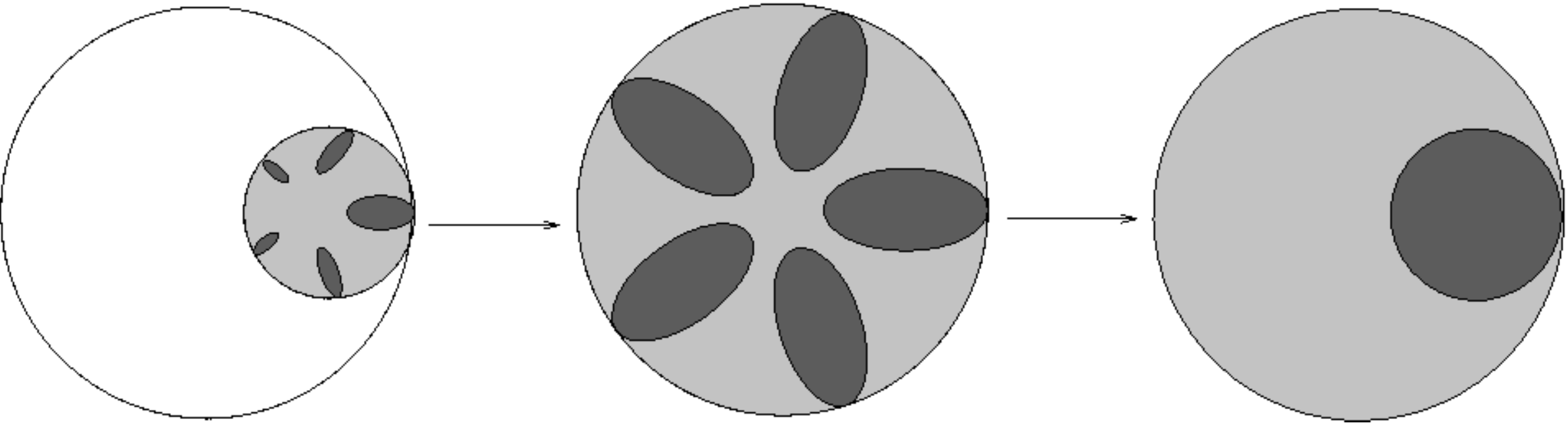}}
\rput(-9.0,1.2){$0$}
\rput(-8.2,1.35){$t_n$}
\rput(-7.5,1.2){$1$}
\rput(-7.1,2.0){$M_n$}
\rput(-5.3,1.3){$0$}
\rput(-3.75,1.2){$1$}
\rput(-3.3,2.0){$z^{q_n}$}
\rput(-1.5,1.2){$0$}
\rput(-0.0,1.2){$1$}
\end{center}
\caption{\it Schematic illustration of the mappings $\ph_n$ on the unit disk as a composition of the M\"obius transformation $M_n$
 and the $q_n^{th}$ power map. The light grey part on the left is sent onto the light grey part on the right, as are the dark gray parts
 (the dark gray disk on the right is close to $0$ so its pre-images are thin).
}
\end{figure}
 This mapping is characterized by the fact that it sends $0$ to $\infty$, $t_n$ to $0$ and $1$ to itself.
We define the sequence of mappings $(\ph_n)_n$ by
$$\ph_n(z)  = \left(M_n(z)\right)^{q_n}.$$

\begin{rem}
\rm{
We will always suppose that $t_n$ belongs to the unit disk.
 Since we suppose $p_n/q_n\tend 0$ this is true for $n$ big enough.}
\end{rem}

 Let
$$\Phi_n = \ph_n\circ\cdots\circ \ph_0.$$
 We denote by $K_\infty$ the set of points of $\barre{\disque}$ which do not escape under $\Phi_n$. That is,
$$K_\infty = \{z\in\barre{\disque}:\forall n,\Phi_n(z)\in\barre{\disque}\}.$$

Note that $$  K_\infty=\bigcap_n K_n,\quad  \text{where } \quad K_n=
\{z\in\barre{\disque}:\,\ph_n\circ\cdots\circ\ph_0(z)\in\barre{\disque}\}. $$
Note that we also have $K_n=\{z\in\barre{\disque}:\forall k=0,\dots,n,\,\ph_k\circ\cdots\circ\ph_0(z)\in\barre{\disque}\}.$
Thus $(K_n)_n$ is a decreasing sequence of non-empty compact sets containing $1$.
In particular, $K_\infty$ itself is non-empty and compact.

The model is defined by the sequence of mappings $(\ph_n)_n$.
 The compact set $K_\infty$ will play an important role in the study of the model and also as part of its realization.

Recall that satellite bifurcations correspond to cycles collisions.
For example, it means that a small perturbation of a polynomial in the
quadratic family having a parabolic fixed point with multiplier
different from $1$ has a cycle which belongs entirely to some small
neighbourhood of the perturbed fixed point.

If the multiplier of the fixed point of the former polynomial is
$e^{2\ii\pi p_n/q_n}$ and if the rotation number of the cycle of the
perturbed polynomial is $p_{n+1}/q_{n+1}$, then the displacement of
the fixed point under perturbation is of the order of magnitude
$p_{n+1}/q_{n+1}$ while the explosion of the cycle happens at a speed
whose order of magnitude is $\left(\frac{p_{n+1}}{q_{n+1}}\right)^{1/q_n}$.

The latter polynomial $f_n=z^2+c_n$, which is a perturbation of a
polynomial with a parabolic fixed point, is renormalizable. The
renormalization replaces the $q_n^{\rm th}$ iterate of the mapping $f_n$ by a
mapping $f_{n+1}=\mathscr{R}f_n$. The new map $f_{n+1}$ has a fixed point which is the image of the exploding
cycle by the renormalization map.

There exists a map $\tilde{\ph}_n$ defined on the domain of
renormalization such that $\mathscr{R}f_n\circ\tilde{\ph}_n=\tilde{\ph}_n\circ f_n^{\circ q_n}$.
In the case where the quadratic polynomial $f_n$ is infinitely
satellite renormalizable, the quadratic polynomial $f_{n+1}$ is again
renormalizable. Thus an infinitely renormalizable $f_0$ yields a
sequence of quadratic maps $(f_n)_n$.

The mapping $\ph_n$ is designed to be a geometric model of the straightening map of the $n^{\rm th}$ satellite renormalization.
In particular it has the following properties :
\begin{itemize}
 \item a fixed point $1$ which represents the critical point ;
 \item the point $t_n$ which represents a fixed point, is sent to the center of the unit disk ;
 \item the set of pre-images of $t_{n+1}$ by $\ph_n$ represents the exploding cycle ;
 \item the power map sends this cycle to a unique fixed point  for the renormalized map.
\end{itemize}

When we consider the quadratic family we can use the Douady-Sullivan criterion (compare \ref{critere de douady-sullivan})
 to show that if the set of accumulation points of the sequence of exploding cycles does not contain the critical point
 then the Julia set of the limit polynomial is not locally-connected.

By its very construction, we know that $K_\infty$ must contain these accumulation points.
This fact allows us to determine a conjectural criterion for non-local connectedness of the Julia set.

\begin{figure}[ht]
\begin{center}
\scalebox{0.12}{\includegraphics{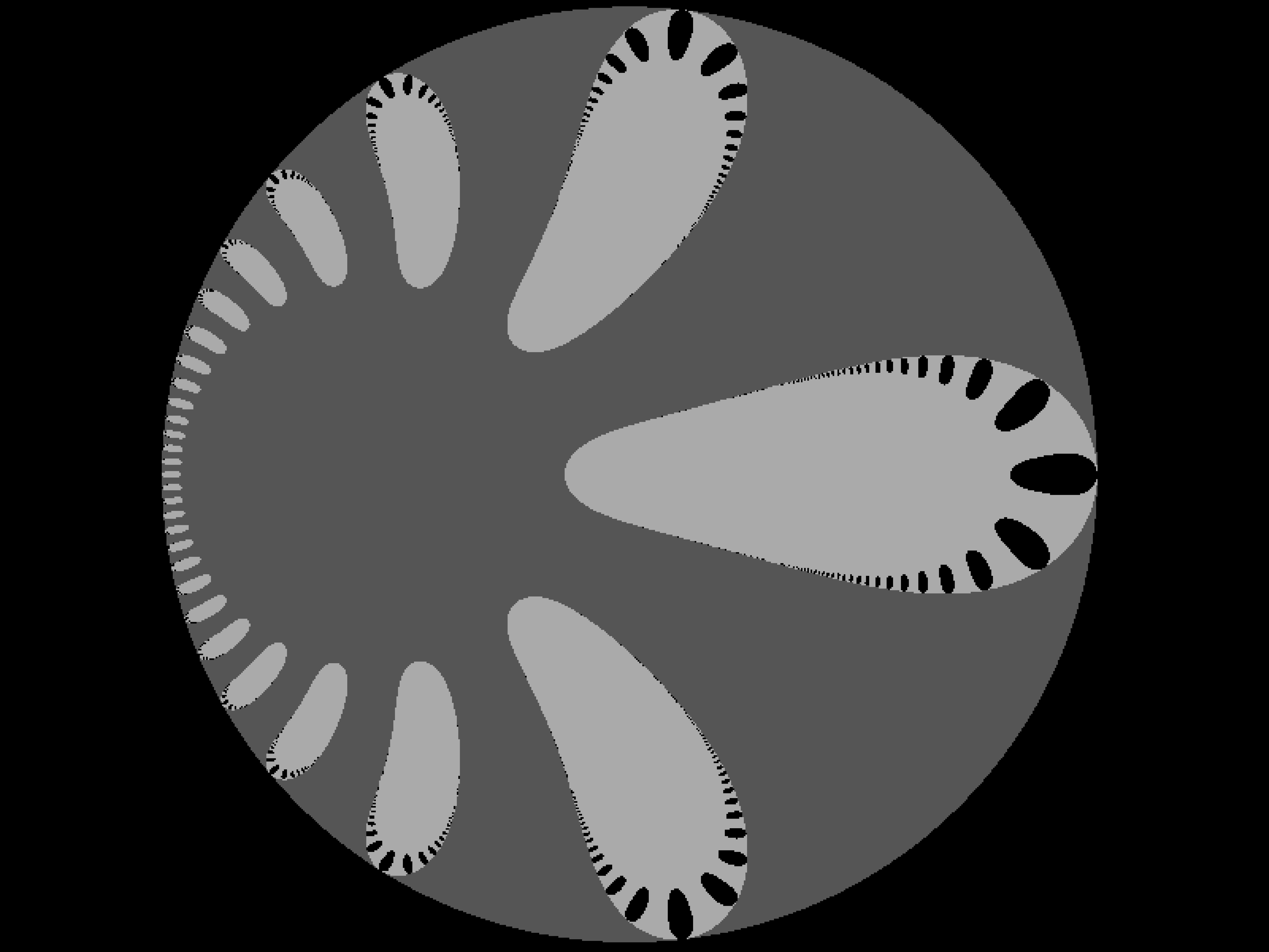}}
\end{center}
\caption{\it Example of the residual set $K_n$ with $n=2$ (it is a magnification of the right part of the unit disk).
Here we have $\fr{p_0}{q_0}=\inv{28}$, $\fr{p_1}{q_1}=\inv{39670}$ (the value of $p_2/q_2$ plays no role in the shape of $K_2$).
The black region surrounding everything is outside any $K_n$, the dark gray represents $K_0\sans K_1$, the light gray $K_1\sans K_2$
and the dark regions inside the light gray is $K_2$.
The compact set $K_n$ is symmetric with respect to the real axis
 and its intersection with the real axis contains a line segment bounded by $1$ on the right (compare section \ref{sec comp critique}).}
\end{figure}

\subsection{The residual compact set}

We begin with the study of the compact set $K_\infty=K_\infty\left((p_n/q_n)_n,C\right)$ called the {\it residual compact set}. Recall that it is defined by
$$ K_\infty = \bigcap_{n=0}^\infty K_n,$$
where $K_n=\{z\in\barre{\disque}:\forall k\leq n,\,\Phi_n(z)\in\barre{\disque}\}$ and $\Phi_n = \ph_n\circ\cdots\circ \ph_0$.

In the next two sections we will label the connected components of $K_\infty$ with an ``odometer'' and prove a result about the topological type of some of its components.

\subsubsection{The address of a point in $K_\infty$}

\begin{lem}\label{lem preimage du disque est un disque a droite}
 Let $t\in]0,1[$ and $M(z)=\fr{1-t/z}{1-t}$. Then $M^{-1}\left(\barre{\disque}\right)$ is the closed disk that has the line segment
 $\left[\fr{t}{2-t},1\right]$ as a diameter.
\end{lem}
\preuve The M\"obius transformation $M$ commutes with $z\mapsto\barre{z}$, sends $1$ to itself and $\fr{t}{2-t}$ to $-1$.\cqfd

\begin{lem}
 Let $n\in\naturels$ and $E_n$ denote the mapping $z\mapsto z^{q_n}$. Then for each connected component of $E_{n}^{-1}\left(M_{n+1}^{-1}\left(\barre{\disque}\right)\right)$
 there exists a unique $k\in\{0,\dots,q_{n}-1\}$ such that $e^{2\pi\ii k\fr{p_{n}}{q_{n}}}$ belongs to this component.
\end{lem}
\preuve From the previous Lemma we know that $M_{n+1}^{-1}\left(\barre{\disque}\right)$ is a disk which lies strictly to the right of $0$.
 Then  the connected components of its preimage by $E_{n}$ are contained in sectors of angles $\fr{\pi}{q_{n}}$ separated by sectors by sectors of the same angle.
Moreover $1\in M_{n+1}^{-1}\left(\barre{\disque}\right)$ so each component contains one and only one $q_{n}^{\rm th}$ root of $1$.
 Finally, note that since $\fr{p_n}{q_n}$ is reduced, the sets $\{e^{2\pi\ii k\fr{p_{n}}{q_{n}}}\mid k\in \N\}$ and $\{e^{2\pi\ii k\fr1{q_{n}}}\mid k\in \N\}$ coincide.  \cqfd

\begin{cor}
 The number of connected components of $K_n$ is $N_n=\ds{\prod_{k=0}^{n-1}}q_k$, with $N_0=1$. 
\end{cor}
\preuve For $n\geq 2$, the mapping $M_{n}$ is a homeomorphism between
 $E_{n-1}^{-1}\left(M_{n}^{-1}\left(\barre{\disque}\right)\right)$ and $\ph_{n-2}\circ\cdots\circ\ph_0(K_n)$,
 and when $n=1$, between $K_1$ and $E_{0}^{-1}\left(M_{1}^{-1}\left(\barre{\disque}\right)\right)$. \cqfd

Thanks to the previous Lemma, we can label the components of
$$\ph_n^{-1}\left(\ph_{n+1}^{-1}\left(\bd\right)\right)=\ph_n^{-1}\left(M_{n+1}^{-1}\left(\bd\right)\right)$$
with the elements of $\relatifs_{q_n}=\relatifs/q_n\relatifs$.
We label the component whose image under $M_n$ contains $e^{2\pi\ii k\fr{p_n}{q_n}}$ by $k\in\relatifs_{q_n }$.

Using this we can define \textit{the address of} $z\in K_n$ by $(k_0,\dots,k_{n-1})\in\relatifs_{q_0}\times\cdots\times\relatifs_{q_{n-1}}$, where the $k_j$ are determined by the condition that 
$\ph_{j-1}\circ\cdots\circ\ph_0(z)$ belongs to the component of $\ph_{j}^{-1}\left(\bd\right)$ which has been labelled $k_{j}$.

The address of a point $z\in K_\infty$ is defined as the infinite sequence $(k_0,\dots,k_n,\dots)\in\ds{\prod_{n=0}^\infty}\relatifs_{q_n}$
 which is such that for each $n$ $(k_0,\dots,k_{n-1})$ is the address of $z$ in $K_n$.
 Every $z\in K_\infty$ has one and only one address but the same address may correspond to several $z$ (see below).

\begin{definition}
 Let $z\in K_\infty$. We say that $\alpha=(\alpha_0,\alpha_1,\dots)\in\ds{\prod_{n\in\naturels}}\relatifs_{q_n}$ is 
\textit{the address of} $z$ \textit{in} $K_\infty$ if for all $n\in\naturels$, the point 
 $M_{n+1}\circ\ph_n\circ\cdots\circ\ph_0(z)$ belongs to the same connected component of $\ph_n^{-1}\left(\bd\right)$ as $e^{2\pi\ii \alpha_n\fr{p_n}{q_n}}$.
\end{definition}

In order to describe the structure of the compact set $K_\infty$ we need to introduce odometers.
Given a sequence of positive integers $N_n$ such that $N_n$ divides $N_{n+1}$, we call the \textit{odometer with scale} $(N_n)_{n\in\naturels}$ the set
$$\mathscr{O}=\relatifs_{N_0}\times\ds{\prod_{n\in\naturels}}\relatifs_{N_{n+1}/N_n}$$
equipped with the product topology of the discrete topology on each factor
 and with a continuous adding map $\sigma:\mathscr{O}\rightarrow\mathscr{O}$
 defined by the following (compare \cite{ref-odometres}) : For all $j=(j_n)_n\in\mathscr{O}$,  $(\sigma(j))_0=j_0+1$ and for $n>0$, 
\rae
\left(\sigma(j)\right)_n = \left\{\begin{array}{l}
j_n+1 \mbox{ if } \forall k\leq n-1,\, j_k=\fr{N_{k+1}}{N_k}-1,\\
j_n \mbox{ otherwise.}
\end{array}\right.
\are
Topologically, an odometer $\mathscr{O}$ is a Cantor set.

In the following we will identify the set of addresses
 $$K_{addr}=\prod_{n\in\naturels}\relatifs_{q_n}$$
 with the odometer with scale $(N_n)_{n\geq 1}$ where $N_n=\ds{\prod_{k=0}^{n-1}}q_k$.
We will refer to this odometer as the addresses odometer.

The adding map might be relevant from the dynamical perspective but not for the study of the topology $K_\infty$;
 compare the proof of Lemma \ref{lem ref pour zq}.

\begin{prop} Let $\mathscr{P}(K_\infty)$ be the set of subsets of $K_\infty$. 
 Let $\pi:K_\infty\rightarrow\mathscr{P}(K_\infty)$ be the mapping that sends a point to the connected component of $K_\infty$ it belongs to. Consider the final  topology on $\mathscr{P}(K_\infty)$  with respect to $\pi$, which is the finest topology that makes the map $\pi$ continuous. \\
Then the set of connected components of $K_\infty$ equipped with the final topology is homeomorphic to the addresses odometer.
\end{prop}
\preuve Let $\alpha\in\ds{\prod_{n\in\naturels}}\relatifs_{q_n}$. Then the set of points which have
 $(\alpha_0,\dots,\alpha_{n-1})$ as their address in $K_n$
 is a connected compact set homeomorphic to the closed unit disk. It  follows that the set of points which have $\alpha$ as their address in $K_\infty$
 is a connected component of $K_\infty$. Hence we have a one-to-one correspondence between the set of connected components of $K_\infty$
 and the set of addresses. We just need to show that the mapping that sends a point to its address is continuous.

Let $z,z'\in K_\infty$ and let $\alpha$, $\alpha'$ be their respective addresses. Suppose that there exists $n$ such that $\forall m\geq n,\,\alpha_m=\alpha'_m$.
Then $\ph_n\circ\cdots\circ\ph_0(z)$ and $\ph_n\circ\cdots\circ\ph_0(z')$ belong to the same connected component of $K_n$.
 This component has a neighbourhood in $\complexes$ which is disjoint from the other components.
 Hence close addresses require the points to be close. \cqfd

\subsubsection{Topology of the critical component}\label{sec comp critique}

In what follows the symbol ``$\arg$'' will denote the argument of a complex number in $]-\pi,\pi]$.

\begin{lem}
 Let $C>1$, $q\in\naturels^*$, $t=C/q$ and let $\ph:\bd\rightarrow\riemann$ be defined by $\ph(z)=\left(\fr{1-t/z}{1-t}\right)^q$.
Then for all $z\in\bd$ such that $|\arg(z-t)|\leq \fr{\pi}{2}$ we have $|\arg\ph(z)|>\fr{C|\parim z|}{2|z|^2}$.
\end{lem}
\preuve Because of the real symmetry of the mapping $\ph$, it is sufficient to show the Lemma for all $z=x+\ii y\in\barre{\disque}$ such that $x\geq t$ and $y>0$.
Under these hypotheses, $\arg(z)=\arcsin\left(y/|z|\right)$ and $\arg(z-t)=\arcsin\left(y/|z-t|\right)$.

Since $\arg\ph(z)=q\left(\arg(z-t)-\arg(z)\right)\geq 0$, we have
$$\arg\ph(z) = q \left(\arcsin\left(\fr{y}{|z-t|}\right)-\arcsin\left(\fr{y}{|z|}\right)\right).$$
The function $\arcsin$ is convex on $|0,1]$, so
\rae
\arcsin\left(\fr{y}{|z-t|}\right) & \geq &
\arcsin\left(\fr{y}{|z|}\right)+\fr{|z|}{x}\left(\fr{y}{|z-t|}-\fr{y}{|z|}\right).
\are
We estimate the difference $\fr{y}{|z-t|}-\fr{y}{|z|}$. Let $r=|z|$. Then
\rae
|z-t|^2 & = & r^2\left(1+\fr{t^2}{r^2}-2\fr{t x}{r^2}\right),
\are
 from which it follows that
\rae
|z-t| & \leq & r\left(1+\demi\left(\fr{t^2}{r^2}-2\fr{t x}{r^2}\right)\right).
\are
Since $\demi\left(\fr{t^2}{r^2}-2\fr{t x}{r^2}\right)\leq  0$, we have
\rae
\inv{|z-t|} & \geq & \inv{r}\fr{1}{1+\demi\left(\fr{t^2}{r^2}-2\fr{t x}{r^2}\right)}\\
 & \geq & \inv{r}\left(1-\demi\left(\fr{t^2}{r^2}-2\fr{t x}{r^2}\right)\right).
\are
As a consequence, $\inv{|z-t|}-\inv{|z|}\geq\fr{xt}{r^3}-\demi\fr{t^2}{r^3}$.

From the above it follows that
\rae
\arg \ph(z)& \geq & \fr{qt y}{r^2}\left(1-\demi\fr{t}{x}\right)\\
 & \geq & \fr{C y}{r^2}\left(1-\demi\fr{t}{x}\right).
\are
But $x\geq t$, so we have $\arg\ph(z)\geq\fr{Cy}{2r^2}$.\cqfd

\begin{cor}\label{cor controle de largument}
 Under the assumptions of the previous Lemma we have
\rae
|\arg(z)|\leq\fr{\pi}{C}|\arg\ph(z)||z|.
\are
\end{cor}
\preuve From the previous Lemma we have $|\arg\ph(z)|\geq\fr{C}{2}\fr{|\parim  y|}{|z|^2}$.
 Note that $\fr{|\parim  y|}{|z|}=\left|\sin(\arg(z))\right|$.
 Using the fact that $|\sin t|\geq\fr{2}{\pi}|t|$ for all $t\in[-\fr{\pi}{2},\fr{\pi}{2}]$, we obtain
 $|\arg\ph(z)|\geq\fr{C}{\pi|z|}|\arg z|$.\cqfd

\begin{definition}\label{def I0}
 Let $K_\infty$ be the residual compact set of the model given by the data $C>1$ and $(p_n/q_n)_n$.
The 
 \textit{ critical component} $I_0$ \textit{of} $K_\infty$ is the connected component of $K_\infty$ which contains $1$.
The 
 \textit{ critical component} $K_{n,0}$ \textit{of} $K_n$ is the connected component of $K_n$ which contains $1$.
\end{definition}
The critical component $I_0$ is the set of points in $K_\infty$ which have $(0,\dots,0,\dots)$ as their address.
The set $K_{n,0}$ is the set of points in $K_n$ which have $(0,\dots,0)$ as their address in $K_n$,
 or, equivalently, whose address in $K_\infty$ starts with $n$ noughts.

\ \\
It follows from the definition that
$$I_0=\bigcap_{n\in\naturels}K_{n,0}.$$

\begin{lem}
 The mapping $M_{n+1}\circ\ph_n\circ\cdots\circ\ph_0:K_{n,0}\rightarrow\bd$
 is the restriction of a biholomorphism defined in a neighbourhood of $K_{n,0}$.
\end{lem}

\preuve 
The mapping $\ph_k$ is the composition of a power map with a M\"obius transformation.
 Then, for all $k<n$, there exists a holomorphic mapping $\psi_k$ defined on $\complexes\sans\reels_-$
 such that $\ph_k\circ\psi_k=\id_{\complexes\sans\reels_-}$.
By definition, for all $k<n$, $\Phi_{k+1}(K_{n,0})\subset\barre{\disque}$, so $\ph_k\circ\cdots\circ\ph_0(K_{n,0})\subset\{\parel >0\}$.
Hence $K_{n,0}\subset\psi_0\circ\cdots\circ\psi_{n-1}\left(\complexes\sans\reels_-\right)$.
\cqfd

\begin{lem}
 The intersection of $K_{n,0}$ with the real axis is a line segment containing $1$.
\end{lem}
\preuve The homeomorphism $\left(M_{n}\circ\ph_{n-1}\circ\dots\circ\ph_0\right)|_{K_{n,0}}$
 is a one-to-one mapping between the real points of $K_{n,0}$ and the real points of $\bd$.\cqfd

\begin{lem}\label{lem restriction sur C}
 Suppose that $C\geq\pi$.
Then for all $n\in\naturels$, 
$$K_{n,0}\subset\left\{z\in\bd:|\arg(z)|\leq\fr{\pi}{2}\left(\fr{\pi}{C}\right)^n\right\}.$$
\end{lem}

Assuming $C\geq\pi$ might not be optimal (compare Lemma \ref{lem inclusions C decroissante}).
We do not know whether this hypothesis is necessary for the above result.

\preuve Let $z\in K_{n,0}$. From Lemma \ref{lem preimage du disque est un disque a droite} it follows that
$\left|\arg(z_k-t_{k+1})\right|\leq \fr{\pi}{2}$ for all $k=0,\dots,n-1$,
 where $z_k=\ph_k\circ\cdots\circ\ph_0(z)\in\barre{\disque}$. 
Then Corollary \ref{cor controle de largument} implies that for all $k\leq n-1$,\linebreak
 $|\arg z_k|\leq\fr{\pi}{C}|\arg z_{k+1}|$. \cqfd

\begin{cor}
 If $C>\pi$, the critical component of $K_\infty$
 is either the point $\{1\}$ or a non trivial line segment on the real axis.
\end{cor}

\begin{rem}\rm{
The components of $K_\infty$ that are sent to $\Phi_n(I_0)$ by some $\Phi_n$ are homeomorphic to $I_0$.
 As a consequence, we showed that a countable dense set of components of $K_\infty$ are homeomorphic to the same line segment 
(possibly reduced to a point).
These are dense because their images by the projection onto the odometer is the dense set of sequences
 $\alpha\in\ds{\prod_{n\in\naturels}}\relatifs_{q_n}$ with finite support ($\alpha_n=0$ for $n$ big enough).}
\end{rem}

The next conjecture illustrates what a topological model for the compact space $K_\infty$ could be.
 Notice that the definition of this model is close to those of a straightbrush and hairy arcs (compare \cite{aarts-oversteegen}, \cite{devaney-straightbrush}).

Recall that $K_{addr}$ is the addresses odometer, homeomorphic to the Cantor set of connected components of $K_\infty$.

\begin{conjecture}
 There exists a closed subset $B$ of $K_{addr} \times [0,1]$ homeomorphic to $K_\infty$ which satisfies
\begin{enumerate}
 \item $K_{addr}\times\{0\}\subset B$;
 \item For all $\alpha\in K_{addr}$ there exists $e_\alpha\in[0,1]$
such that
 $(\alpha,t)\in B$ if and only if $0\leq t\leq e_\alpha$.
(Such $e_\alpha$ is called an upper endpoint);
 \item The set of upper endpoints is dense in $B$.
\end{enumerate}

\end{conjecture}

The definition of $K_\infty$ might recall in some points one of the definitions of a Cantor Bouquet
 (see the characterization of the Julia set of some maps of the exponential family as a Cantor Bouquet in \cite{devaney-straightbrush}) :
\begin{itemize}
\item The set $K_\infty$ is the set of non-escaping points under the compositions of a ordered countable family of holomorphic mappings $(\ph_n)_n$.
\item It is an intersection of a decreasing sequence of compact sets $K_n$ (in the Riemann sphere).
\item For every connected component $\omega$ of $K_n$ the cardinality of the set of connected components of $K_{n+1}$ which are included in $\omega$ are the same.
\item The mapping $\ph_n$ sends $K_{n+1}\cap\omega$ homeomorphically onto $K_n$ for all connected component $\omega$ of $K_n$.
\item Connected components are ordered (vertical lines for the Julia sets of the exponential maps, circular ordering for the present object of our study).
 It allows to label the connected components of $K_\infty$.
\end{itemize}

The above conjecture shows how to unroll and how to straighten $K_\infty$.

\subsection{The conjectural non-local-connectedness criterion}

\begin{lem}\label{lem x0 in limit set}
Let
$$x_0=\inf\{x\in[0,1]:\forall y>x\,\forall n\in\naturels,\,\Phi_n(y)\geq t_{n+1}\}.$$

\ \\
Then,
\begin{itemize}
 \item $[x_0,1]\subset K_\infty$;
 \item $x_0=\ds{\lim_{n\tend\infty}}s_n,$ where, for all $n\in\naturels$,
$s_n$ is the unique pre-image of $0$ by $\Phi_n$ which belongs to the component $K_{n,0}$;
 \item For all $\eps>0$ there exists 
 $y\in [x_0-\eps,x_0[\sans K_\infty$.
\end{itemize}
\end{lem}

\preuve
From Lemma \ref{lem preimage du disque est un disque a droite}
 it follows that if a point $x\in K_n$ is such that $\parel\Phi_n(x)<0$ then $x\nin K_\infty$.
Since the mapping $\ph_n$ is increasing on $[t_n,1]$ and $\Phi_{n-1}(s_n)=t_n$, the mapping $\Phi_n$ is also increasing on $[s_n,1]$.
As a consequence $[x_0,1]\subset K_\infty$ and $x_0$ is the limit of the increasing sequence $(s_n)_n$.

Let $\eps>0$ and let $n$ be such that $x_0-s_n\leq\demi{\eps}$.
If we take $y\in K_{n,0}\cap [s_n-\eps/2,s_n[$, then we have $\Phi_n(y)<0$ so 
$y\in[x_0-\eps,x_0[\sans K_\infty$.
\cqfd

\begin{definition}
 Let $C>1$. We denote by $\mathscr{C}_C$ the set of all parameters $c\in\complexes$
 such that the quadratic polynomial $z^2+c$ is infinitely satellite renormalizable with the sequence of rotation numbers
 $(p_n/q_n)_n$ satisfying the following :
\begin{enumerate}
 \item The sequence of positive numbers $(t_n)_n$ defined by $t_n=C\left|p_n/q_n\right|$ is such that
  $\forall n\in\naturels^*$, $t_n\in]0,1[$.
 \item Let $\ph_n$ be the mapping defined by $\ph_n(z)=\left(\fr{1-t_n/z}{1-t_n}\right)^{q_n}$ and let
 $\Phi_n= \ph_n\circ\cdots\circ\ph_0$. Then
\begin{equation}\label{eq cond phix etc}
\exists x_0\in]0,1[,\forall n\in\naturels,\Phi_{n}(x_0)\geq t_{n+1}.
\end{equation}
\end{enumerate}

\end{definition}

Because of the monotonicity of the mappings $\ph_n$ on $[t_n,1]$,
 the second condition implies that the line segment $[x_0,1]$ is included in $K_\infty$ (compare Lemma \ref{lem x0 in limit set}).
Conversely if there is $x_0\in]0,1[$ such that $[x_0,1]\subset K_\infty$ then $x_0$
 satisfies the condition in (\ref{eq cond phix etc}).

\begin{lem}\label{lem inclusions C decroissante}
 Let $C\geq C'>1$. Then $\mathscr{C}_C\subset\mathscr{C}_{C'}$.
\end{lem}

\preuve Let $c\in\mathscr{C}_C$ and let $\gamma=C'/C$. By the assumption the sequence of real numbers 
 $t_n=C|p_n|/q_n$ is such that
 $t_n\in]0,1[$. Also there exists  $x_0\in]0,1[$ such that for all $n$,
 $\Phi_n(x_0)\geq t_{n+1}$.
 We define the sequence of real numbers $x_n$ by $x_{n+1}=\ph_n(x_n)$.

 Let $t'_n=\gamma t_n\in]0,1[$,
  $\ph'_n(z)=\left(\fr{1-t'_n/z}{1-t'_n}\right)^{q_n}$ and $\Phi'_n=\ph'_n\circ\cdots\circ\ph'_0$.
 We show that the sequence $(x'_n)_n$ defined by $x'_0=x_0$ and $x'_{n+1}=\ph'_n(x'_n)$ satisfies $x'_n\geq x_n$ for all $n$.

By induction, suppose $x'_n\geq x_n$ and $x'_n\in[0,1]$. Since $t'_n\leq t_n$ we have
$$x'_{n+1}\geq \ph_n(x'_n) \geq \ph_n(x_n)=x_{n+1}.$$
Hence $\Phi_n'(x_0)\geq \Phi_n(x_0)\geq t_{n+1}\geq t'_{n+1}$.\cqfd

\vspace{3mm}

The main conjecture is the following (we purposely state it only in the case where $p_n=1$).

\begin{conjecture}[(non-local-connectedness criterion)]\label{conj isr jnlc}
 Let $c\in\complexes$ be such that the quadratic polynomial $Q_c(z)=z^2+c$ is infinitely satellite renormalizable with the sequence of rotation numbers $(1/q_n)_n$.
Then the Julia set of $Q_c$ is not locally-connected if and only if  there exists $C>1$ and a renormalization $Q_{c*}$ of $Q_c$ such that $c*\in\mathscr{C}_C$.
\end{conjecture}

We can justify this conjecture by the Douady-Sullivan criterion \ref{critere de douady-sullivan}.

The cycles are modeled by the centers of the components of $K_n$,
 that is by the preimages of $0$ under $\Phi_n$ for some $n$.
 Thanks to Lemma \ref{lem x0 in limit set} we know that the left boundary of the segment $I_0$ (the critical component) is a limit point of the set of centers of components of $K_n$, $n\in\naturels$.
The other limit points belong to other components which are at a positive distance from the point $1$.
As a consequence, if this segment is not reduced to a point
 the critical point is at a positive distance of the limit set of the cycles and we can apply the Douady-Sullivan criterion.

\ \\

Another conjecture related to this model is the existence of an invariant compact set inside the Julia set which is homeomorphic to $K_\infty$.

\begin{conjecture}\label{conj any 418}
 Let $c\in\complexes$ be such that the quadratic polynomial $Q_c(z)=z^2+c$ is infinitely satellite renormalizable with the sequence of rotation numbers $(1/q_n)_n$.
 
 Then there exists $C>1$, a renormalization $Q_{c*}$ of $Q_c$ and an invariant compact subset of the Julia set of $Q_{c_*}$ homeomorphic to the residual compact set $K_\infty$
  of the model associated to the data $\left(C,(p_n/q_n)_n\right)$.

Moreover, the componentwise dynamics in this set is given by the adding map of the addresses odometer.
\end{conjecture}

Non-locally-connected quadratic Julia sets are still not well understood.
Proving this conjecture may provide valuable information on the structure of and the dynamics on  the Julia set in the case it is not locally-connected.

It would be interesting to know if this homeomorphism
 extends to the whole plane and if it is even quasi-conformal.

\subsubsection{A test of the conjectural criterion}

The following is an example of a situation where the Julia set is not locally-connected and the conditions of the conjecture \ref{conj any 418} are satisfied
 (that is $c_*\in\mathscr{C}_C$ for some renormalization $Q_{c_*}$).

As we mentioned earlier, the article \cite{So} of
 S{\o}rensen 
 does not contain explicit conditions for non-local-connectedness, but
 Milnor has proposed such a condition in \cite{MiNlc}.
G. Levin found an explicit condition which implies Milnor's condition. Indeed Levin's criterion is more general :

\begin{thm}[(Levin, \cite{Levin-2009})]\label{thm de levin}
 Let $c\in\complexes$ be such that the quadratic polynomial $Q_c(z)=z^2+c$ is infinitely satellite renormalizable with the sequence of rotation numbers $(p_n/q_n)_n$. Suppose that
$$\limsup_{n\tend\infty}\left|\fr{p_{n+1}}{q_{n+1}}\right|^{1/q_n}<1.$$
Then the Julia set of $Q_c$ is not locally-connected.
\end{thm}

Indeed the work of Levin yields a more general condition which is not easy to work with.
 No other explicit criterion based on the rotation numbers is known yet.

\begin{thm}
 Let $c\in\complexes$ be such that the quadratic polynomial $Q_c(z)=z^2+c$ is infinitely satellite renormalizable with the sequence of rotation numbers $(p_n/q_n)_n$. Suppose that
\begin{itemize}
 \item the sequence $(p_n)_n$ is bounded while $q_n\tend\infty$,
 \item $\ds{\limsup_{n\tend\infty}}\left|\fr{p_{n+1}}{q_{n+1}}\right|^{1/q_n}<1$.
\end{itemize}
Then for all $C>1$, there exists a renormalization $Q_{c_*}$ of $Q_c$ such that $c_*\in\mathscr{C}_C$.
\end{thm}

\preuve Since we can renormalize, we may suppose that there is an $\alpha\in[0,1[$ such that \linebreak
 $|p_{n+1}/q_{n+1}|\leq\alpha^{q_n}$.
By hypothesis $\fr{p_n}{q_n}\tend 0$, thus we can renormalize so that for all $n$, $t_n<1$, where $t_n$ is defined as $t_n=C|p_{n}/q_{n}|$.

Let $\beta\in]1,1/\alpha[$ and $\eta=\inv{1-\beta\alpha}$.
For the same reason as above, we may also suppose that for all $n$, $\left(\fr{\beta}{1-t_n}\right)^{q_n}\geq C\eta$.

Define the sequence $(x_n)_n$ in the following way. We set $x_0=\eta t_0$ and for all $n\geq 0$ define\linebreak
  $x_{n+1}=\left(\fr{1-t_n/x_n}{1-t_n}\right)^{q_n}$.
Then the sequence $(x_n)_n$ satisfies $x_n\geq \eta t_n$ for all $n$. In fact by induction,
\rae
x_{n+1} & \geq & \left(\fr{1-1/\eta}{1-t_n}\right)^{q_n}=\left(\fr{\beta}{1-t_n}\right)^{q_n}\alpha^{q_n}\\
 & \geq & \eta t_{n+1}
\are
\cqfd

\end{document}